\documentclass[5p, twocolumn]{elsarticle}

\usepackage{multirow}
\usepackage{eurosym}
\usepackage{lscape}
\usepackage{makeidx}
\usepackage{epsfig}
\usepackage[cmex10]{amsmath}
\usepackage{array}
\usepackage{mdwmath}
\usepackage{mdwtab}
\usepackage{eqparbox}
\usepackage[font=footnotesize]{subfig}
\usepackage{url}
\usepackage{multicol}
\usepackage[T1]{fontenc}
\usepackage[charter]{mathdesign}
\usepackage{eulervm}
\usepackage{enumitem}
\usepackage{xcolor}
\usepackage{algorithm,algpseudocode}
\biboptions{square,numbers,sort&compress}
\usepackage{color,soul}
\definecolor{lightgray}{gray}{0.6}
\sethlcolor{lightgray}

\usepackage{graphicx}
\usepackage{caption}
\usepackage{epstopdf}

\journal{International Journal of Electrical Power \& Energy Systems}

\begin{document}

\begin{frontmatter}


\title{Hierarchical Interactive Architecture based on Coalition Formation for Neighborhood System Decision Making}

\author[Virginia Tech]{Jaber Valinejad}
\ead{Jabervalinejad@vt.edu
Corresponding author}
\author[Newcastle]{Mousa Marzband}
\author[LLNL]{Mert Korkali }
\author[Virginia Tech]{Yijun Xu}
\author[UAE]{Ameena Saad Al-Sumaiti}

\address[Virginia Tech]{Bradley Department of Electrical and Computer Engineering, Virginia Polytechnic Institute and State University, Northern Virginia Center, Falls Church, VA 22043, USA}

\address[Newcastle]{Faculty of Engineering and Environment, Department of Physics and Electrical Engineering, Northumbria University Newcastle, Newcastle upon Tyne NE1 8ST, UK}
\address[LLNL]{Computational Engineering Division, Lawrence Livermore National Laboratory, Livermore, CA 94550 USA.
}
\address[UAE]{Advanced Power and Energy Center, Electrical Engineering and Computer Science, Khalifa University, Abu Dhabi, UAE}


\begin{abstract}
In recent years, significant efforts for improving the technical and the economic performance of smart grids have been applied with the presence of different players making decisions in these grids. 
This paper proposes a bilevel hierarchical structure for designing and planning distributed energy resources in home microgrids (H-MGs). 
Such a small-scale grid inside the market environment permits energy exchange among distributed energy resources and home microgrids through a pool market. In this paper, each home microgrid's planner maximizes the performance of distributed energy resources while exchanging energy with other H-MGs. The problem is formulated as a high-level problem decomposed into a set of low-level market clearing problems. 
The global optimal performance in terms of energy cost is met for a market structure (H-MGs, consumers, and retailers) at a Nash optimum point for the formulated scheduling game considering local and general constraints on the spot market. In general, the upper-level structure is based on H-MG generation competition for maximizing their individual and/or group income in the process of forming coalition with other H-MGs. In addition, the functionality of the lower-level of the hierarchical structure is governed by a market clearing based on the price response method by all the DERs enabling H-MGs to change the spot market strategic behavior. Prices are obtained as dual variables of power balance equations. 
This paper will investigate a set of bilevel games, including nonlinear programming problems solved through a complementary problem method. 
Using the binary theory, the bilevel hierarchical structure will be replanned using a nonlinear problem. 
Results prove that the proposed structure can increase the players's profit.

\end{abstract}



\end{frontmatter}





\section*{Nomenclature}
\vspace{-0.2cm}
\begin{tabular}{p{1.7cm} p{6.5cm}}

\multicolumn{2}{c}{$\textbf{Acronyms}$} \\

CHP	& combined heat and power  \\

EB	& electrical boiler \\

EHP	& electrical heat pump \\
ES+/ES- & energy storage (ES) during charging/discharging mode \\
GB	& gas boiler \\

k+/k- & buying/ selling power from/to H-MG \textit{i}/ Retailer $k$ \\
MCP	& market clearing price \\



DR+, DR-	& amount of responsive load demand (RLD) that goes from/comes to other time period to/from $t$ \\

STP	& solar thermal panel \\

TES	& thermal energy storage \\

\end{tabular}

\begin{tabular}{p{1.7cm} p{6.5cm}}




WT	& wind turbine \\
RET & retailer\\
A, B, C & home microgrids \\

\multicolumn{2}{c}{$\textbf{Indices}$} \\

w  & the number of scenario ($w\in\{1,2,\dots,W\}$) \\
$i/m$  & the number of H-MG ($i\in\{1,2,\dots,I\}$)\\
$j$  & the number of DERs existing in H-MG i ($j\in\{1,2,\dots,J\}$) \\
$t$  & time ($t\in\{1,2,\dots,T\}$)\\
$D$ &	the amount of electrical load consumer DERs's demand (i.e., EB, EHP, ES) and/or consumers existing in H-MG $i$ \\
k  & the number of a retailer ($k\in\{1,2,\dots,n\}$)\\
e/h	& electrical/thermal \\
FU	& fuel consumed by thermal DERs (i.e., CHP, EB) \\

\end{tabular}

\begin{tabular}{p{1.7cm} p{6.5cm}}

\multicolumn{2}{c}{$\textbf{Constants}$} \\

$\zeta_{\text{ij}}^{\text{CHP,h}}$ & the thermal efficiency of CHP  \\

$\zeta_{\text{ij}}^{\text{CHP,e}}$/$\zeta_{\text{ij}}^{\prime \text{CHP,e}}$ & the electrical efficiency of CHP  \\
 
$\pi^{\text{FU}}$ & the offer price of natural gas \\
 
$\text{SOC}_{\text{ij}}^{\text{ES,ini}}$ $\text{SOC}_{\text{ij}}^{\text{ES,end}}$ & SOC initial and final value related to ES $j$ in H-MG $i$ [kW] \\
 
$\text{SOC}_{\text{ij}}^{\text{TES,ini}}$ $\text{SOC}_{\text{ij}}^{\text{TES,end}}$ & SOC initial and final value related to TES $j$ in H-MG $i$ [kW] \\


$\zeta_{\text{ij}}^{\text{EB}}$, $\zeta_{\text{ij}}^{\text{TES}}$ & the thermal efficiency of EB/TES \\

$\text{COP}^{\text{EHP}}$ & coefficient of performance of EHP \\

$\overline{\text{MC}}_{\text{ij}}$ & marginal cost of DER $j$ in H-MG $i$ [\pounds/kWh] \\
$\underline{P}_{\text{ij}}$, $\overline{P}_{\text{ij}}$ & maximum electrical/thermal power generated by DER $j$ in H-MG $i$ [kW] \\

\multicolumn{2}{c}{$\textbf{Parameters}$} \\

$\tilde{\lambda_{\text{tw}}}^{\text{MCP}}$ & the value of predicted electrical clearing price at time $t$ under scenario $w$ [\pounds/kWh] \\

$\tilde{P}_{\text{itw}}^{\text{D,e}}$ & predicted load power consumed by the consumers in H-MG $i$ at time $t$ under scenario $w$ [kW] \\

\multicolumn{2}{c}{$\textbf{Functions}$} \\

$\mathbb{R}_{ij}$ & the revenue of DER $j$ related to H-MG $i$ [\pounds] \\

$\mathbb{R}^{i}$ & the revenue of H-MG $i$ in independent operation [\pounds] \\

$\mathbb{R}^{im}$ & the revenue resulting from coalitional operation between H-MG $i$ and H-MG $m$ [\pounds] \\

\multicolumn{2}{c}{$\textbf{Decision Variables}$} \\

$\pi_{\text{ijtw}}^{X,e}$/$\pi_{\text{ijtw}}^{Y,h}$  & electrical/ thermal selling price bid by DER $j$ under scenario $w$ at time $t$ [\pounds/kWh] \\


$P_{\text{ijtw}}^{e}$/$P_{\text{ijtw}}^{h}$ & electrical/thermal power generated by DER $j$ at time $t$ under scenario $w$ [kW] \\


$P_{\text{iktw}}^{e}$ & electrical power sold to/bought from H-MG $i$/ Retailer $k$ at time $t$ under scenario $w$ [kW] \\

$P_{\text{itw}}^{D,e}$ &  active power consumed by the consumers at H-MG $i$ at time $t$ under scenario $w$ [kW] \\

$P_{\text{itt'w}}^{\prime D,e}$ & shifted load power from time interval $t$ to time interval $t'$ at H-MG $i$ under scenario $w$ [kW] \\

$P_{\text{ijtw}}^{\text{EB},e}$, $P_{\text{ijtw}}^{\text{EHP},e}$ & consumed load power by EB/EHP $j$ at H-MG $i$ at time $t$ under scenario $w$ [kW] \\

$\lambda_{t}^{\text{MCP},e}$, $\lambda_{t}^{\text{MCP},h}$ & electrical and thermal market clearing price at time $t$ [\pounds/kWh] \\

$P_{\text{imtw}}^{e-}$, $P_{\text{imtw}}^{e+}$ & buying/selling electrical power from H-MG $i$ to H-MG $m$ [kW] \\

$P_{\text{iktw}}^{e-/+}$/$P_{\text{kitw}}^{e-/+}$ & buying/ selling electrical power from/to H-MG $i$/ Retailer $k$ [kW] \\

$P_{\text{imtw}}^{h-}$/$P_{\text{imtw}}^{h+}$ & buying/selling thermal power from/to H-MG $i$/H-MG $m$ [kW] \\

$\pi_{\text{iktw}}^{e-}$/$\pi_{\text{iktw}}^{e+}$ & supply bids for buying/ selling electrical power from/to H-MG $i$/Retailer $k$ [kW] \\

$\pi_{\text{imtw}}^{h-}$/$\pi_{\text{imtw}}^{h+}$ & supply bids for buying/selling thermal power from/to H-MG $i$/H-MG $m$ [kW] \\



\end{tabular}


\section{Introduction}
\label{sec:Introduction}

Demand-side management (DSM) topics are focused on energy consumption control at the consumer side~\cite{Kohn2015,Liu2015_1}. Such energy control is coordinated by electric utilities, companies and enterprises without controlling distributed energy resources \cite{Li2017,Zazo2017}. When the latter is controlled, the topic would be defined as an energy management \cite{Ma2014}. \par
demandThis objective of this paper is to first propose a base framework for the demand of consumers encompassing H-MGs, and second investigate profits that can be made from operating H-MGs independently or in a coalitional structure in a daily electricity market \cite{Pei2016,Zhang2017}. For these purposes, energy exchange among H-MGs is formulated as a scheduling game. In this paper, competitive monopolies are modeled and simulated in a formulated convex optimization problem \cite{Liu2017,Liu2017_1}. These monopolies are based on three contradictory objectives: H-MGs's and retailers's income maximization, consumers's cost reduction, and demand peak reduction \cite{Panwar2017}. \par
With the aim of a convex optimization problem, a bilevel hierarchical interactive architecture (BL-HIA) algorithm on the condition of reaching a maximum profit is proposed for both consumer and the power generator side \cite{Kinhekar2016,Basu2012}. The optimum performance problem is presented for all DERs existing in multiple H-MGs as a BL-HIA so that it is a mixed integer programming problem \cite{Lee2014,Jayaweera2015}. The upper-level targets maximizing H-MGs's profit through energy exchange among H-MGs as well as H-MGs and retailers for a central optimal performance of a decision maker \cite{Jayaweera2015,Wang2012}. On the other hand, the lower-level of the hierarchal structure of the problem represents an equilibrium problem incorporating DSM for an optimal performance of multiple H-MGs \cite{Kohn2015,Ma2014}. This way, a central optimum performance decision maker for energy optimum exchange among H-MGs with each other and with retailers in an independent and a coalitional performance and with the aim of reaching H-MGs's maximum profit is included at an upper-level decision maker while considering an independent or a coalitional performance of all H-MGs \cite{Giusti2014,Ni2016,Gabbar2016}. The interaction between the two levels of the hierarchal structure of the game is a factor of searching for the optimal solution at both levels \cite{Zazo2017,Ma2014}. Considering the optimal scheduling of all  H-MGs and DERs existing in them in Multiple H-MGs requires solving mathematical program with equilibrium constraints (MPEC) equivalent to bilevel problem. This bilevel problem can be looked at as a multiple-leader-common-follower game. The aim of implementing this game structure is finding a final equilibrium point in which none of H-MGs and consumers can increase their profit by changing in the generation and the consumption schedules. Furthermore, the BL-HIA structure accounts for decisions resulting from forming a coalition among H-MGs to maximize the profit and also exchange energy among them. \par
The contribution in this work can be summarized as follows: First, the proposed BL-HIA structure is preferable over the proposed structure in \cite{Marzband2015ETEP} as it is a multi-ownership structure that permits, forming coalition among H-MGs and explicitly increasing the competition among H-MGs and consumers rather than an independent operation of H-MGs. Second, the BL-HIA is adequate for modelling problems with several leaders (i.e., H-MGs) having their own individual objective functions when operating independently or in a coalitional manner (upper-level problem). Such a game is to optimize several followers (i.e., consumers inserted in the bilevel structure). These models are related to situations where actions and followers's performance in BL-HIA have a significant effect on decisions made by leaders. This fact is related to the case in which H-MGs's profit (a leader) depends on the amount of energy which is sold for supplying consumers existing in the power grid (as follower). The general view of the hierarchical structure and optimization problems has been shown with the proposed model in Figure.~\ref{fig:Fig_1}.
Third, a better strategy maximizing consumers's satisfaction in terms of demand supply and H-MGs's profit is presented in comparison to a single level structure. Finally, the BL-HIA structure is solved by formulating an equivalent one level mixed-integer nonlinear programming (MINLP) problem deploying the KKT optimization conditions. The innovations in this paper can be summarized as follows: \par

\begin{enumerate}[noitemsep,nolistsep]
    \item An optimum programming solution within H-MGs generation as a BL-HIA structure; 
    \item Providing a multiple-leader-common-follower game which states the effectiveness of the market competition in multiple H-MGs through solving a BL-HIA structure;
    \item Developing a new model for demand-side management;
    \item Accommodating both DR resources and storage devices in the market operation to achieve a comprehensive solution exploiting all flexibilities; and
    \item Proposing an advanced electricity market for active distribution networks based on game theory;
\end{enumerate}

\begin{figure}
\centering
\includegraphics[width= \columnwidth]{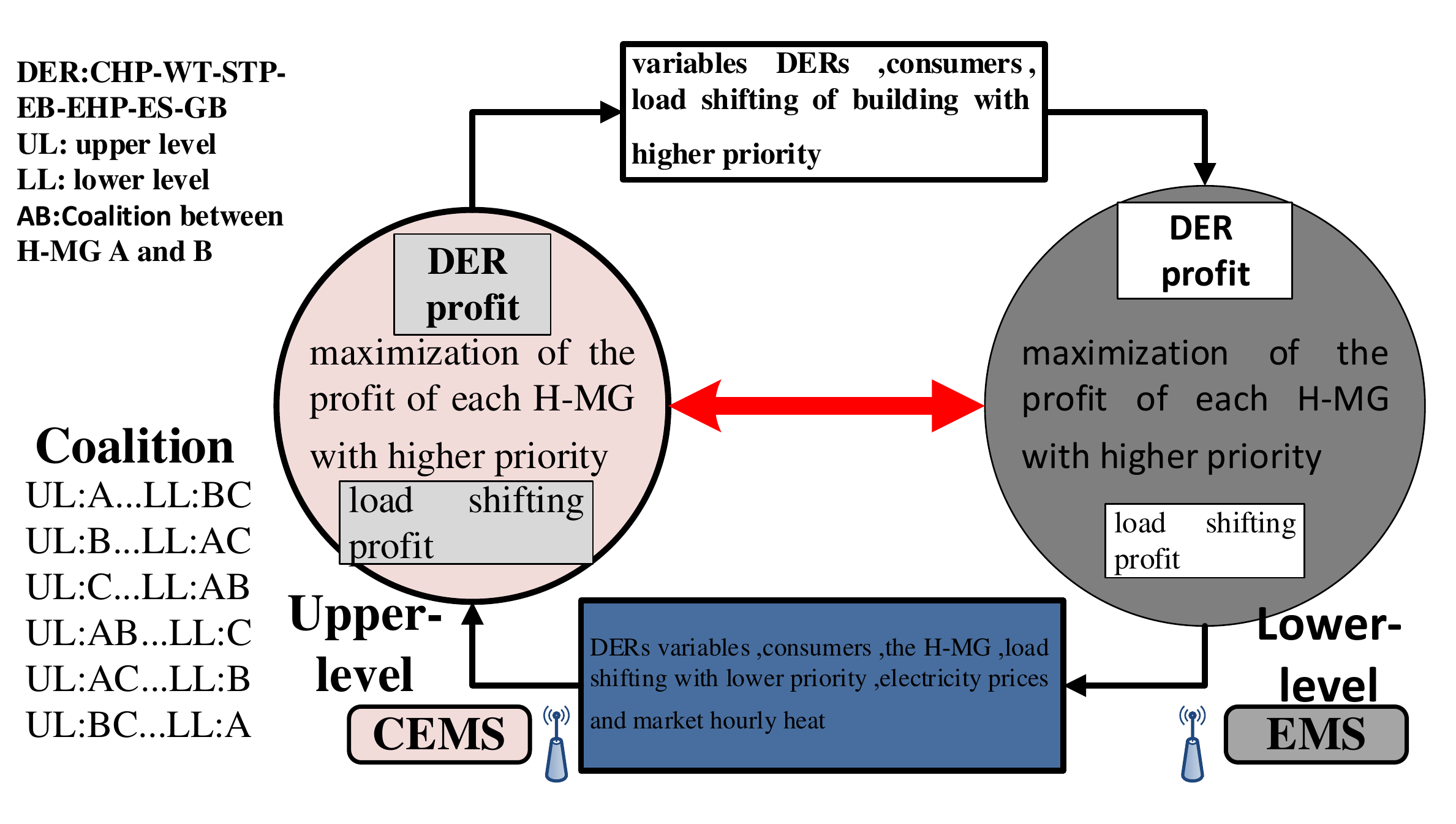}
\caption{The proposed BL-HIA structure illustrating a variety of coalition formation among H-MGs. UL:A...LL:BC implies that the HMG A is modeled in the upper level while HMGs B and C are modeled in lower level.Distributed energy resources encompass set of distributed generation, electrical and thermal energy storage.}
\label{fig:Fig_3}
\end{figure}

\section{The Proposed BL-HIA Structure }
\label{sec:bilevel model hierarchical structure}

\begin{figure}
  \centering
  \includegraphics[width= \columnwidth]{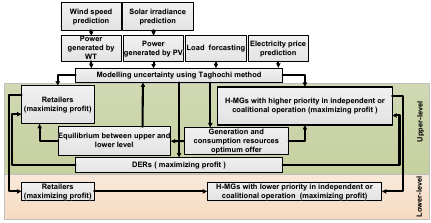}
  \caption{BL-HIA structure.}
  \label{fig:Fig_1}
  \end{figure}

The problem encountered by the H-MGs for an independent or a coalitional operation can be modelled as a bilevel structure that is a decision-making problem including several agents which try to optimize their corresponding objective functions on a connectable dependent set. An agent is, in fact, an object which can act as a DER or connected to other units. The BL-HIA structure is shown in Figure.~\ref{fig:Fig_1}.
The upper-level problem states maximizing the profit of H-MGs having a higher priority of operating in an independent or a coalitional operation on the condition of: 1) satisfying upper-level constraints and 2) satisfying a set of lower-level problems. H-MGs with a higher priority of operation in the upper-level problem are identified based on their price bids. Upper-level constraints include limits on the quantity and supply bids of DERs resources, the minimum accessible power capacity by the market regulator, and buying/selling quantities by H-MGs and retailers. Lower-level problem states the market clearing prices (MCP) with the aim of maximizing the profit of H-MGs having lower priority of operation subject to meeting  equilibrium constraints for each H-MG, generation/consumption limits, and the number of consumers participating in the DSM program. \par
As it is observed in Figure.~\ref{fig:Fig_1}, the higher priority of H-MGs operation (in an independent or a coalitional operation) is defined by the energy excess/shortage gap and supply/demand bids that permits a maximization of an expected profit. Maximizing the profit of each H-MG is achieved by considering the fact that each agent at the lower-level problem shows an optimal operation in correspondence to the income of H-MGs with higher priority of operation based on the offered price. This optimal operation includes an estimation of the demand supplied and shifted by each consumer in an independent or a coalitional operation. It must be emphasized that the competition among H-MGs with higher and lower priority (or competitor H-MGs) are explicitly modelled at upper and lower-levels. It must be noted that upper-level and lower-level problems shown in Figure.~\ref{fig:Fig_1} become related to each other. In other words, lower-level problems estimate the price and the quantity of competitor H-MGs which directly affect the profit of H-MGs inserted in the above problem. In other words, decisions related to forming a coalition  and taking a bidding strategy by the competitor strategic H-MGs in the upper-level problem have also a significant effect on the MCP resulted from the lower-level problem.\par

The BL-HIA structure is shown in Figure.~\ref{fig:Fig_3}. In BL-HIA structure is also considered for modelling the uncertainty of pool prices, electrical/thermal load demand, and H-MGs buying/selling prices. The proposed bilevel model has been simplified as a one-level problem using KKT method \cite{valinejad2017generation} for the sake of convexification. 


 
\section{Decision-Making Process by BL-HIA Structure}

\begin{figure}
 \centering
 \includegraphics[width= \columnwidth]{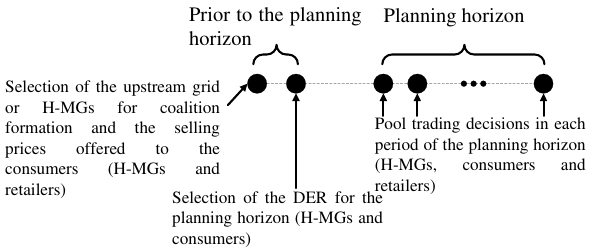}
 \caption{Decision-making process.}
 \label{fig:Fig_2}
 \end{figure}

The decision-making process in the BL-HIA structure of H-MGs, consumers and retailers can be summarized as shown in Figure.~\ref{fig:Fig_2}. At the beginning of the scheduling horizon, each H-MG presents the necessary decisions on DSM for an independent or a coalitional operation with other H-MGs. Moreover, supply/demand bids are provided to the consumers during this horizon. These decisions are made under uncertainties of a future pool market prices, consumers's load profile, and competitor H-MGs's supply bids. The supply bids are a function of the cost of DER installed in the H-MG.

\begin{enumerate}[noitemsep,nolistsep]
\item Consumers's choice of energy provider:\par
 When each H-MG offers a suuply bid, consumers are to choose a H-MG as an energy provider to their  electrical/thermal load  during the scheduling horizon. These decisions are made based on reliable information on such prices estimated under uncertainties of pool prices  and demand. For modelling purposes, several sets of consumers are created by grouping consumers with similar specifications responding to H-MG's offered prices).
\item Energy exchange in a pool market by the H-MGs:\par
 After stabilizing the H-MGs's performance (an independent or a coalitional operation) and setting supply and demand bids, each H-MG can decide in each time interval of the scheduling horizon on the quantity (to/from other H-MGs in the pool market) to supply the demand of their consumers.

\end{enumerate}

 \section{Problem Formulation}

 
The H-MGs's scheduling  problem is formulated in a BL-HIA structure. It must be noted that dual variables have been separated by comma after equality and inequality constraints. This section will briefly present  models deployed for load shifting, and those representing the interaction among DER, H-MGs, consumers and retailers as well as the coalition among H-MGs. Then, BL-HIA problem formulation will be presented.

 \subsection{DR Objective Function and Constraints}

\noindent
\begin{equation}
\label{eq:1}
\textbf{max} \;\; 
\mathbb{R}_{i}^{\text{DSM}}=
0.5
\times
\displaystyle{\sum_{t'=1}^{T}}
P_{\text{it'w}}^{D,e}
\times
\left(
\tilde{\lambda}_{\text{tw}}^{\text{MCP}}
- \tilde{\lambda}_{\text{t'w}}^{\text{MCP}}
\right)
\end{equation}

The profit resulting from the participation of consumers in DSM program is calculated from~(\ref{eq:1}):

 \noindent
 \begin{equation}
 \label{eq:2}
 P_{\text{itw}}^{D,e}
 =
 \tilde{P}_{\text{itw}}^{D,e}
 +
 \displaystyle{\sum_{t'}}
 P_{\text{itt'w}}^{\prime D,e}:
 \;\;
 \gamma_{\text{itw}}^{D,e}
 \end{equation}

 \noindent
  \begin{equation}
  \label{eq:3}
  \displaystyle{\sum_{t \to t'}}
  P_{\text{itt'w}}^{\prime D,e}
  \leq
  \tilde{P}_{\text{itw}}^{D,e}:
  \;\;
  \overline{\eta}_{\text{itw}}^{\text{D,e}}
  \;\;
  \forall(\tilde{\lambda}_{\text{tw}}^{\text{MCP}}> \tilde{\lambda}_{\text{t'w}}^{\text{MCP}}), \forall P_{\text{tt'w}}^{\prime D,e}> 0
  \end{equation}

  \noindent
    \begin{equation}
    \label{eq:4}
  P_{\text{itt'w}}^{\prime D,e}> 0: \underline{\eta}_{\text{itt'w}}^{D,e}
    \end{equation}

\noindent
    \begin{equation}
    \label{eq:5}
  P_{\text{itt'w}}^{\prime D,e}= 0: \eta_{\text{itt'w}}^{D,e},\forall t=t'
    \end{equation}

Equation~(\ref{eq:3}) states that if the value of the demand shifted from time interval $t$ to $t'$ is given when $\tilde{\lambda}_{\text{tw}}^{\text{MCP}}> \tilde{\lambda}_{\text{t'w}}^{\text{MCP}}$, then this value is not to exceed the predicted load value at $t$ (i.e., $\tilde{P}_{\text{itw}}^{D,e}$). Equations~(\ref{eq:4}) and (\ref{eq:5}) are set to ensure that load shifting is not defined for the same interval. The load shifted from $t$ to $t'$ is equivalent to the negative value of the demand deducted from the predicted demand at $t$ as described by

  \noindent
      \begin{equation}
      \label{eq:6}
    P_{\text{itt'w}}^{\prime D,e}= -P_{\text{itt'w}}^{\prime D,e}: \eta_{\text{itt'w}}^{\prime D,e}.
      \end{equation}

\subsection{CHP Objective Function and Constraints}

\noindent
      \begin{equation}
      \label{eq:7}
\textbf{max} \;\; 
\mathbb{R}_{ij}^{\text{CHP}}=    
  \displaystyle{\sum_{t=1}^{T}}
  P_{\text{ijtw}}^{\text{CHP,e}}
  \times
  \pi_{\text{ijtw}}^{\text{CHP,e}}
  + P_{\text{ijtw}}^{\text{CHP,h}}
    \times
    \pi_{\text{ijtw}}^{\text{CHP,h}}
    - \frac{P_{\text{ijtw}}^{\text{CHP,h}}}{N_{i}^{\text{CHP,h}}}
    \times
    \pi^{\text{FU}}
      \end{equation}

The objective is to maximize the profit that can be made through CHP systems's participation in DSM program as in~(\ref{eq:7}):

\noindent
\begin{equation}
\label{eq:8}
\underline{P}_{\text{ij}}^{\text{CHP,e}}
\leq
P_{\text{ijtw}}^{\text{CHP,e}}
\leq
\overline{P}_{\text{ij}}^{\text{CHP,e}}:
\;\;
\underline{\eta}_{\text{ijtw}}^{\text{CHP,e}},
\overline{\eta}_{\text{ijtw}}^{\text{CHP,e}},
\end{equation}

\noindent
\begin{equation}
\label{eq:9}
\underline{P}_{\text{ij}}^{\text{CHP,h}}
\leq
P_{\text{ijtw}}^{\text{CHP,h}}
\leq
\overline{P}_{\text{ij}}^{\text{CHP,h}}:
\;\;
\underline{\eta}_{\text{ijtw}}^{\text{CHP,h}},
\overline{\eta}_{\text{ijtw}}^{\text{CHP,h}},
\end{equation}

\noindent
\begin{equation}
\label{eq:10}
P_{\text{ijtw}}^{\text{CHP,e}}
= \zeta_{\text{ij}}^{\text{CHP,e}}
\times
\frac{P_{\text{ijtw}}^{\text{CHP,h}}}{\zeta_{\text{ij}}^{\text{CHP,h}}}
+
\zeta_{\text{ij}}^{\prime \text{CHP,e}}:
\;\;
\gamma_{\text{ijtw}}^{\text{CHP,e}}
\end{equation}

Equations~(\ref{eq:8}) and (\ref{eq:9}) state upper and lower limits on the power generated by the CHPs. Equation~(\ref{eq:10}) describes the power generated by the CHPs as a function of the system's efficiency

\subsection{WT Objective Function and Constraints}

\noindent
\begin{equation}
\label{eq:13}
\textbf{max} \;\; 
\mathbb{R}_{ij}^{\text{WT}}=    
\displaystyle{\sum_{t=1}^{T}}
P_{\text{ijtw}}^{\text{WT}}
\times
\pi_{\text{ijtw}}^{\text{WT}}
\end{equation}

The profit resulting from the participation of WT in the DSM program is calculated by~(\ref{eq:13}):

\noindent
\begin{equation}
\label{eq:14}
0
\leq
P_{\text{ijtw}}^{\text{WT}}
\leq
\overline{P}_{\text{ij}}^{\text{WT}}:
\;\;
\underline{\zeta}_{\text{tw}}^{\text{WT,e}},
\overline{\zeta}_{\text{tw}}^{\text{WT,e}}
\end{equation}

Equation~(\ref{eq:14}) (and similar constraints for determining the limit on DER resources) state programmed power generation of DER controllable and non-controllable resources. Equation~(\ref{eq:14}) is related to the wind turbine electrical power whose maximum limit is a parameter having uncertainty.

\subsection{ES and TES Objective Functions and Constraints}

\noindent
\begin{equation}
\label{eq:16}
\textbf{max} \;\; 
\mathbb{R}_{ij}^{\text{ES/TES}}=    
\displaystyle{\sum_{t=1}^{T}}
P_{\text{ijtw}}^{\text{ES/TES}}
\times
\pi_{\text{ijtw}}^{\text{ES/TES}}
\end{equation} 

The profit made by the participation of an electrical/thermal energy storage (ES/TES)  system in a DSM program is measured by the objective function in~(\ref{eq:16}):

\noindent
\begin{equation}
\label{eq:17}
-\underline{P}_{\text{ij}}^{\text{ES/TES}}
\leq
P_{\text{ijtw}}^{\text{ES/TES}}
\leq
\overline{P}_{\text{ij}}^{\text{ES/TES}}:
\;\;
\underline{\eta}_{\text{tw}}^{\text{ES/TES}}, \overline{\eta}_{\text{tw}}^{\text{ES/TES}}
\end{equation}

\noindent
\begin{equation}
\label{eq:18}
\text{SOC}_{\text{ijtw}}^{\text{ES/TES}}
= \text{SOC}_{\text{ij(t-1)w}}^{\text{ES/TES}}
+ \frac{P_{\text{ijtw}}^{\text{ES/TES}}}{\overline{P}_{\text{ij}}^{\text{ES/TES}}}:
\;\;
\gamma_{\text{ijtw}}^{\text{ES/TES}}
\end{equation}

\noindent
\begin{equation}
\label{eq:19}
\underline{\text{SOC}}_{\text{ij}}^{\text{ES/TES}}
\leq
\text{SOC}_{\text{ijtw}}^{\text{ES/TES}}
\leq
\overline{\text{SOC}}_{\text{ij}}^{\text{ES/TES}}:
\;\;
\underline{\eta}_{\text{ijtw}}^{\text{ES/TES,SOC}}, \overline{\eta}_{\text{ijtw}}^{\text{ES/TES,SOC}}
\end{equation}

\noindent
\begin{equation}
\label{eq:20}
\text{SOC}_{\text{ij(t-1)w}}^{\text{ES/TES}}
= \text{SOC}_{\text{ij}}^{\text{ES/TES,ini}}:
\;\;
\gamma_{\text{ijtw}}^{\prime \text{ES/TES}}
\end{equation}

\noindent
\begin{equation}
\label{eq:21}
\text{SOC}_{\text{ij(t=24)w}}^{\text{ES/TES}}
= \text{SOC}_{\text{ij}}^{\text{ES/TES,end}}:
\;\;
\gamma_{\text{ijw}}^{\prime \prime \text{ES/TES}}
\end{equation}

Equations~(\ref{eq:17})--(\ref{eq:21}) are the equations governing the operation of ES/TES systems. The operation of ES/TES system is subject to generation limits as in~(\ref{eq:17}) and state of charge limits as in~(\ref{eq:18})--(\ref{eq:21}). It should be noted that~(\ref{eq:19}) states the charge/discharge rate of the ES/TES system.

\subsection{EB Constraints}

\noindent
\begin{equation}
\label{eq:23}
\textbf{max} \;\; 
\mathbb{R}_{ij}^{\text{EB}}=    
\displaystyle{\sum_{t=1}^{T}}
\left(
P_{\text{ijtw}}^{\text{EB,h}}
\times
\pi_{\text{ijtw}}^{\text{EB,h}}
- P_{\text{ijtw}}^{\text{EB,e}}
\times
\pi_{\text{ijtw}}^{\text{EB,e}}
\right)
\end{equation}

The profit made by the participation of EB in the DSM program is calculated by~(\ref{eq:23}):

\noindent
\begin{equation}
\label{eq:24}
P_{\text{ijtw}}^{\text{EB,h}}
= \zeta_{\text{ij}}^{\text{EB}}
\times
P_{\text{ijtw}}^{\text{EB,e}}:
\;\;
\zeta_{\text{tw}}^{\text{EB,h}}
\end{equation}

\noindent
\begin{equation}
\label{eq:25}
0
\leq
P_{\text{ijtw}}^{\text{EB,h}}
\leq
\overline{P}_{\text{ij}}^{\text{EB,h}}:
\;\;
\underline{\eta}_{\text{ijtw}}^{\text{EB,h}}, \overline{\eta}_{\text{ijtw}}^{\text{EB,h}}
\end{equation}

Eqs.~(\ref{eq:24}) and (\ref{eq:25}) state the consumed amount of electrical power and the generated heat in the electrical boiler, respectively.

\subsection{EHP Objective Function and Constraints}

\noindent
\begin{equation}
\label{eq:27}
\textbf{max} \;\; 
\mathbb{R}_{ij}^{\text{EHP}}=    
\displaystyle{\sum_{t=1}^{T}}
\left(
P_{\text{ijtw}}^{\text{EHP,h}}
\times
\pi_{\text{ijtw}}^{\text{EHP,h}}
- P_{\text{ijtw}}^{\text{EHP,e}}
\times
\pi_{\text{ijtw}}^{\text{EHP,e}}
\right)
\end{equation}

The profit made by the participation of EHP in the DSM program is calculated by~(\ref{eq:27}):

\noindent
\begin{equation}
\label{eq:28}
P_{\text{ijtw}}^{\text{EHP,h}}
= \text{COP}^{\text{EHP}}
\times
P_{\text{ijtw}}^{\text{EHP,e}}:
\;\;
\gamma_{\text{ijtw}}^{\text{HP,h}}
\end{equation}

\noindent
\begin{equation}
\label{eq:29}
0
\leq
P_{\text{ijtw}}^{\text{EHP,h}}
\leq
\overline{P}_{\text{ijtw}}^{\text{EHP,h}}:
\;\;
\underline{\eta}_{\text{ijtw}}^{\text{EHP,h}}, 
\overline{\eta}_{\text{ijtw}}^{\text{EHP,h}}
\end{equation}

Equation~(\ref{eq:28}) describes the relation between the consumed amount of electrical and thermal power generated by the thermal pump. In addition,~(\ref{eq:29}) represents limits on the thermal generation.

\subsection{GB Objective Function and Constraints}

\noindent
\begin{equation}
\label{eq:38}
\textbf{max} \;\; 
\mathbb{R}_{ij}^{\text{GB}}=    
\displaystyle{\sum_{t=1}^{T}}
\left(
P_{\text{ijtw}}^{\text{GB}}
\times
\pi_{\text{ijtw}}^{\text{GB}}
- \frac{P_{\text{ijtw}}^{\text{GB}}}{N_{\text{ij}}^{\text{GB}}}
\times
\pi^{\text{FU}}
\right)
\end{equation}

The profit made by the participation of GB in the DSM program is calculated by~(\ref{eq:38}).

\noindent
\begin{equation}
\label{eq:39}
0
\leq
P_{\text{ijtw}}^{\text{GB}}
\leq
\overline{P}_{\text{ij}}^{\text{GB}}:
\;\;
\underline{\eta}_{\text{ijtw}}^{\text{GB,h}},
\overline{\eta}_{\text{ijtw}}^{\text{GB,h}}
\end{equation}

Equation~(\ref{eq:39}) states the allowable limits on the heat generation by the gas boiler.

\subsection{STP Objective Function and Constraints}

\noindent
\begin{equation}
\label{eq:41}
\textbf{max} \;\; 
\mathbb{R}_{ij}^{\text{STP}}=    
\displaystyle{\sum_{t=1}^{T}}
\left(
P_{\text{ijtw}}^{\text{STP}}
\times
\pi_{\text{ijtw}}^{\text{STP}}
\right)
\end{equation}

The profit made by the participation of STP in a DSM program is calculated by Eq.~(\ref{eq:41}):

\noindent
\begin{equation}
\label{eq:420}
0
\leq
P_{\text{ijtw}}^{\text{STP}}
\leq
\overline{P}_{\text{ij}}^{\text{STP}}:
\;\;
\underline{\eta}_{\text{ijtw}}^{\text{STP,h}},
\overline{\eta}_{\text{ijtw}}^{\text{STP,h}},
\end{equation}


Equation~(\ref{eq:420}) states the heat generation allowable limits for the operation of a thermal solar panel. In common with a wind turbine, the solar panel thermal maximum limit is also considered to be an uncertainty factor.

\subsection{DERs's Price Bid Constraints}

\noindent
\begin{equation}
\label{eq:11}
0
\leq
\pi_{\text{ijtw}}^{\text{X,e}}
\leq
\tilde{\lambda}_{\text{tw}}^{\text{MCP,e}}
\;\;
\end{equation}

\noindent
\begin{equation}
\label{eq:12}
0
\leq
\pi_{\text{ijtw}}^{\text{Y,h}}
\leq
2\times\tilde{\lambda}_{\text{tw}}^{\text{MCP,e}}
\end{equation}   

Equations~(\ref{eq:11}) and (\ref{eq:12}) are related to the electrical and thermal price bids governing DERs's operation where $X$ includes CHP, ES, WT; and $Y$ encompasses CHP, EB, EHP, TES, GB, STP. However, the value of the upper and lower bounds can vary with respect to the system deployed.

\section{H-MGs's Independent and Coalitional Operation}

Two scenarios are implemented to simulate the performance of the proposed BL-HIA structure. These scenarios are described as follows:

\begin{itemize}
\item Scenario 1:

\noindent
This scenario describes independent operations of H-MGs (\{A\},\{B\},\{C\},\{RET\}). A single-level algorithm is deployed to model this scenario as further clarified by the independent operation of H-MGs.

\item Scenario 2: 
(\{A,\{B,C\}\}, \{\{A,B\},C\}, \{\{A,C\},B\}, \{B,\{A,C\}\}, \{\{B,C\},A\}, \{C,\{A,B\}\}):

\noindent

This scenario describes a coalition among H-MGs taking place at a single level of the BL-HIA structure and operating in an independent operation at the other level. A representation of such a scenario can take the shape of (\{A,BC\}, \{AB,C\}, \{AC,B\}, \{B,AC\}, \{BC,A\}, \{C,\{AB\}\}).  This scenario also investigates the effect of  the  lower-level H-MGs forming a coalition on changing the strategy of operating the upper-level H-MG with a high priority independently. The mathematical model of this scenario is further clarified by the coalitional operation of H-MGs as shown next.

\end{itemize}

\subsection{Scenario 1}

Equation~(\ref{eq:44}) states the profit gained by selling energy by Retailer $k$ to all H-MGs.

\noindent
\begin{equation}
\label{eq:44}
\textbf{max} \;\; 
\mathbb{R}_{k}^{\text{RET}}=    
\displaystyle{\sum_{t=1}^{T}}
\displaystyle{\sum_{i=1}^{n}}
\left(
P_{\text{kitw}}^{\text{e-}}
\times
\pi_{\text{kitw}}^{\text{e-}}
- P_{\text{iktw}}^{\text{e-}}
\times
\pi_{\text{iktw}}^{\text{e-}}
\right)
\end{equation}

Equation~(\ref{eq:45}) state the profit obtained through the independent operation of the H-MGs (\{A\}, \{B\},\{C\}). It is worth mentioning that if  a DER does not exist in an H-MG, it is not considered in the respective objective function.

\noindent
\begin{equation}
\label{eq:45}
\begin{array}{l}
\mathbb{R}^{\text{A/B/C}}=    
\displaystyle{\sum_{t=1}^{T}}
\displaystyle{\sum_{w=1}^{W}}
\displaystyle{\sum_{j=1}^{J}}
\left(
\begin{array}{l}
P_{\text{ijtw}}^{\text{ES}}
\times
\pi_{\text{ijtw}}^{\text{ES}}
+ P_{\text{ijtw}}^{\text{WT}}
\times
\pi_{\text{ijtw}}^{\text{WT}} \\
P_{\text{ijtw}}^{\text{STP}}
\times
\pi_{\text{ijtw}}^{\text{STP}}
+ P_{\text{ijtw}}^{\text{EHP,h}}
\times
\pi_{\text{ijtw}}^{\text{EHP,h}} \\
P_{\text{ijtw}}^{\text{TES}}
\times
\pi_{\text{ijtw}}^{\text{TES}}
+ P_{\text{ijtw}}^{\text{EB,h}}
\times
\pi_{\text{ijtw}}^{\text{EB,h}} \\
- P_{\text{ijtw}}^{\text{EB,e}}
\times
\pi_{\text{ijtw}}^{\text{EB,e}} 
+ P_{\text{ijtw}}^{\text{GB}}
\times
\pi_{\text{ijtw}}^{\text{GB}}  \\
+ P_{\text{ijtw}}^{\text{CHP,e}}
\times
\pi_{\text{ijtw}}^{\text{CHP,e}} 
+ P_{\text{ijtw}}^{\text{CHP,h}}
\times
\pi_{\text{ijtw}}^{\text{CHP,h}} \\
+ P_{\text{imtw}}^{\text{e-}}
\times
\pi_{\text{imtw}}^{\text{e-}} 
+ P_{\text{iktw}}^{\text{e-}}
\times
\pi_{\text{iktw}}^{\text{e-}} \\
+ P_{\text{imtw}}^{\text{h-}} 
\times
\pi_{\text{imtw}}^{\text{h-}} 
- P_{\text{imtw}}^{\text{e+}}
\times
\pi_{\text{imtw}}^{\text{e+}} \\
+ P_{\text{iktw}}^{\text{e+}}  
\times
\pi_{\text{iktw}}^{\text{e+}} 
- P_{\text{imtw}}^{\text{h+}}
\times
\pi_{\text{imtw}}^{\text{h+}}\\
\end{array}
\right) \\
+ \mathbb{R}_{i}^{\text{DSM}} 
\;\; (i=A/B/C)
\end{array}
\end{equation}

\subsection{Scenario 2}

Equations~(\ref{eq:48})--(\ref{eq:510}) state the profit obtained through the coalitional operation of H-MGs (\{A\},\{B\},\{C\}) at an upper-level or a lower-level.

\noindent
\begin{equation}
\label{eq:48}
\textbf{max} \;\; 
\mathbb{R}_{i}^{\text{\{AB\}}}=
\mathbb{R}_{i}^{\text{\{A\}}}
@
\mathbb{R}_{i}^{\text{\{B\}}}
\end{equation}



\noindent
\begin{equation}
\label{eq:510}
\textbf{max} \;\; 
\mathbb{R}_{i}^{\text{\{ABC\}}}=
\mathbb{R}_{i}^{\text{\{A\}}}
@
\mathbb{R}_{i}^{\text{\{B\}}}
@
\mathbb{R}_{i}^{\text{\{C\}}}
\end{equation}


\noindent where @ states the coalition among different H-MGs. In addition, coalitional scenario of (\{B,AC\}) means that the first part (B) is related to the objective function defined at the upper-level and a second part (AC) is  related to coalition between H-MGs A and C but defined at a lower level.

Equations~(\ref{eq:51})--(\ref{eq:56}) state the profit obtained from the coalitional operation of H-MGs A, B, and C at an upper or a lower level. \par
The right-hand side of these relations is made up of two parts. The first part is related to the objective function defined at the upper level; whereas, the second part is related to the objective function defined at the lower level.

\noindent
\begin{equation}
\label{eq:51}
\textbf{max} \;\; 
\mathbb{R}_{i}^{\text{\{A,BC\}}}=
\mathbb{R}_{i}^{\text{\{A\}}}
@
\mathbb{R}_{i}^{\text{\{BC\}}}
\end{equation}

\noindent
\begin{equation}
\label{eq:52}
\textbf{max} \;\; 
\mathbb{R}_{i}^{\text{\{AB,C\}}}=
\mathbb{R}_{i}^{\text{\{AB\}}}
@
\mathbb{R}_{i}^{\text{\{C\}}}
\end{equation}

\noindent
\begin{equation}
\label{eq:53}
\textbf{max} \;\; 
\mathbb{R}_{i}^{\text{\{AC,B\}}}=
\mathbb{R}_{i}^{\text{\{AC\}}}
@
\mathbb{R}_{i}^{\text{\{B\}}}
\end{equation}

\noindent
\begin{equation}
\label{eq:54}
\textbf{max} \;\; 
\mathbb{R}_{i}^{\text{\{B,AC\}}}=
\mathbb{R}_{i}^{\text{\{B\}}}
@
\mathbb{R}_{i}^{\text{\{AC\}}}
\end{equation}

\noindent
\begin{equation}
\label{eq:55}
\textbf{max} \;\; 
\mathbb{R}_{i}^{\text{\{BC,A\}}}=
\mathbb{R}_{i}^{\text{\{BC\}}}
@
\mathbb{R}_{i}^{\text{\{A\}}}
\end{equation}

\noindent
\begin{equation}
\label{eq:56}
\textbf{max} \;\; 
\mathbb{R}_{i}^{\text{\{C,AB\}}}=
\mathbb{R}_{i}^{\text{\{C\}}}
@
\mathbb{R}_{i}^{\text{\{AB\}}}
\end{equation}

\section{Mathematical Formulation of BL-HIA Structure}

In the upper-level problem, each H-MG seeks to maximize its amount of profit. Each upper-level problem's objective function states the income of each H-MG with a lower-level priority for different scenarios. These objective functions which must be maximized have been defined as the sum of the product of electrical/thermal price offers and electrical/thermal powers sold to consumers of each H-MG minus the cost of operation of DERs. \par
BL-HIA structure includes upper-level problems and a set of lower-level problems under each scenario w. It is noted that if a DER is considered in the upper-level, their constraints from~(\ref{eq:1})--(\ref{eq:12}) are considered in the upper-level. similarly  is the case for lower-level. 
The upper-level problem includes decision making regarding the possibility of forming coalition among H-MGs and their supply bids to achieve a higher profit. On the other hand, the quantity by DERs/consumers resources along DSM program are included in the lower-level problem. It must be noted that all the power exchange among H-MGs and retailers's decisions are to be made on the upper-level problem. In comparison, decision-making variables at the lower-level include all the power generation by DER resources. The upper-level objective function is considered after maximizing the income of retailer or H-MGs in the case of independent or in the coalitional operations with other H-MGs (under investigated scenarios). \par
The income of H-MGs is defined as the product of the proposed offer for selling power to H-MGs [\pounds/kWh] by the amount of power sold to them [kWh] minus the product of a power purchase offer from H-MGs [kWh] by the amount of power bought from quantity provided by other H-MGs.

\subsection{The Upper-Level Problem}


The upper-level relationship is formulated in this section. The formula expressing the DER relationship given in~(\ref{eq:1})--(\ref{eq:12}) are applicable if the related DER is considered to be in the upper level.


\subsubsection{Objective function}
As was stated earlier, the objective functions of the upper-level and the lower-level problems can be in the form of~(\ref{eq:45})--(\ref{eq:510}). Here, the profit obtained from the coalition or the independent operation of H-MGs with more priority describes the objective function of the upper-level problem. The upper-level problem is to maximize the expected profit to be made by  each H-MG in the case of an individual or a group operation as well as a retailer. 




\subsubsection{Upper-level problem constraints}

\noindent
\begin{equation}
\label{eq:57}
(\ref{eq:1})-(\ref{eq:12})
\end{equation}

Equations~(\ref{eq:1})--(\ref{eq:12}) are for H-MGs with higher priority. It is very obvious if an H-MG with a higher priority does not include any of the mentioned DERs, then the corresponding constraints of such DERs are to be excluded from the problem formulation. 
Each H-MG in the upper-level problem makes strategic decisions as per the following:

\begin{itemize}
\item Decisions of DERs on supply bid in the lower level of the problem;
\item Price strategic offering decisions of consumers in a H-MG on price offers.
\end{itemize}


\subsection{The Lower-Level Problem}

Each lower-level problem states maximizing the profit of an individual H-MG or a group of H-MGs with a lower priority over different scenarios. The objective of the lower-level problem is to increase the profit of DERs. Thus, CEMS 
is to reduce the operating cost given limitations ruling over each of the players (i.e., H-MGs, retailer, consumers). Players in BL-HIA structure declare the amount of their generated power and supply bids offered to the CEMS. After simulating the bilevel problem, the electrical and thermal energy prices, the amount of the quantity by each of the players are provided.

\subsubsection{Objective function}


Considering the cases described in~(\ref{eq:45})--(\ref{eq:510}), the objective function of the lower-level problem can be taken from the numerator of the lower-level objective function.


\subsubsection{Lower-level problem constraints load shifting and the DERs}

Equations~(\ref{eq:1})--(\ref{eq:12}) apply for each H-MG with a lower priority. It is very obvious that if an H-MG with a lower priority for example has no CHP, then its constraints must not be considered.

\subsection{Power Exchange Constraints between H-MGs and Retailers}

\noindent
\begin{equation}
\label{eq:60}
0
\leq
P_{\text{iktw}}^{e-}
\leq
\eta_{i}
\times
P_{\text{itw}}^{e-}:
\;\;
\overline{\eta}_{\text{iktw}}^{\text{Grid,e}}, \underline{\eta}_{\text{iktw}}^{\text{Grid,e}}
\end{equation}

\noindent
\begin{equation}
\label{eq:61}
0
\leq
P_{\text{iktw}}^{e+}
\leq
\eta_{i}
\times
P_{\text{itw}}^{e+}:
\;\;
\overline{\eta}_{\text{iktw}}^{\text{Grid,e}}, \underline{\eta}_{\text{iktw}}^{\text{Grid,e}}
\end{equation}



Equations~(\ref{eq:60}) and (\ref{eq:61}) show the allowable limits on power exchange between retailers and H-MGs. 

\subsection{Electrical Balance Constraint}

\noindent
\begin{equation}
\label{eq:62}
\begin{array}{l}
P_{\text{iktw}}^{e-}+ P_{\text{itw}}^{D,e}
-
\left(
P_{\text{ijtw}}^{\text{CHP,e}}
+ P_{\text{ijtw}}^{\text{WT}}
+ P_{\text{ijtw}}^{\text{ES}}
- P_{\text{ijtw}}^{\text{EB,e}}
- P_{\text{ijtw}}^{\text{EHP,e}}
+ P_{\text{iktw}}^{e+}
\right) \\
+ \displaystyle{\sum_{m=1}^{n}} P_{\text{imtw}}^{e+}
- \displaystyle{\sum_{m=1}^{n}} P_{\text{mitw}}^{e+}
= 0: \;\; \lambda_{\text{tw}}^{e}
\end{array}
\end{equation}

Equation~(\ref{eq:62}) states the equilibrium relation between H-MGs generated and consumed electrical power and electrical power exchange (KW) with retailers. The MCP in the grid is equal to the dual variable of~(\ref{eq:62}).



\subsection{Thermal Balance Constraints}

\noindent
\begin{equation}
\label{eq:63}
\begin{array}{l}
P_{\text{itw}}^{D,h}
- P_{\text{ijtw}}^{\text{CHP,h}}
+ P_{\text{ijtw}}^{\text{EB,h}}
+ P_{\text{ijtw}}^{\text{EHP,h}}
+ P_{\text{ijtw}}^{\text{TES}}
+ P_{\text{ijtw}}^{\text{GB}}
+ P_{\text{ijtw}}^{\text{STP}} \\
+ \displaystyle{\sum_{m=1}^{n}} P_{\text{imtw}}^{h-}
- \displaystyle{\sum_{m=1}^{n}} P_{\text{imtw}}^{h+}
= 0: \;\; \lambda_{\text{tw}}^{h}
\end{array}
\end{equation}

Equation~(\ref{eq:63}) states the relation between the thermal generated and consumed power. The thermal power price in the dual variable grid corresponds to~(\ref{eq:63}). \par
After the determination of price offers related to electricity and heat and also the amount of electrical and thermal power generation and consumption of each player, the profit made by each of the players is determined.
















Since each one of these lower-level problems is continuous and convex, it may be shown by its specific constraints including Karush-Kuhn-Tucker (KKT) conditions \cite{Valinejad2017}. 
By using KKT conditions, the constraints for an independent or a coalition operation of H-MGs include the following cases:

\begin{itemize}
\item Primal constraints~(\ref{eq:1})--(\ref{eq:12});
\item Equality constraints obtained from the derivative of a Lagrange expression relative to lower-level variables 
\item Complementary constraints obtained based on lower-level inequalities~(\ref{eq:3}),(\ref{eq:4})--(\ref{eq:8})--(\ref{eq:9})--(\ref{eq:14})--(\ref{eq:17})--(\ref{eq:19})--(\ref{eq:25})--(\ref{eq:29})- (\ref{eq:39})--(\ref{eq:420})--(\ref{eq:11})--(\ref{eq:12})--(\ref{eq:60})--(\ref{eq:61}) 

\end{itemize}







\section{Results and Discussion}

\begin{figure}
   \centering
   \includegraphics[width= \columnwidth]{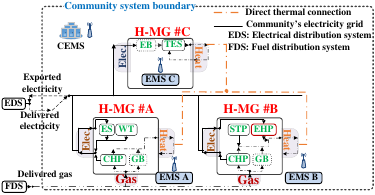}
   \caption{The grid under study.}
   \label{fig:Fig_4}
   \end{figure}

The grid under study is shown in Figure.~\ref{fig:Fig_4}. The energy storage systems installed in H-MGs ({A} and {C}) are for storing excess electrical and thermal energy generation. The capacity and the number of installed equipment in each H-MG are shown in Table~\ref{tab:1}.  
 
\begin{table}
\centering
\caption{The Capacity and the Number of Equipment Installed in Each H-MG}
\begin{tabular}{|l| c| c| c|}
\hline
DERs&	H-MG A &	H-MG B &	H-MG C \\
\hline
CHP	& & &\\	
\hline	
Electrical output ($\text{kW}_{\text{ele}}$) &	142	& 207	& - \\
\hline	
Thermal output ($\text{kW}_{\text{th}}$) &	104 &	140	& - \\
\hline
EHP ($\text{kW}_{\text{th}}$) &	- &	700	 & -  \\
\hline
WT ($\text{kW}_{\text{ele}}$) &	50	& -	& - \\
\hline
STP ($\text{kW}_{\text{ele}}$) &	- &	600	 & - \\
\hline
ES ($\text{kW}_{\text{ele}}$) (2kWh) &	- &	500	 & - \\
\hline
TES ($m^{3}$)	&	- & - &	4 \\
\hline
GB ($\text{kW}_{\text{th}}$) &	2$\times$150 &	2$\times$150 &	- \\
\hline
EB ($\text{kW}_{\text{th}}$) &	-  & - &	2$\times$100  \\
\hline
\end{tabular}
\label{tab:1}
\end{table}

\begin{figure}[h!] 
    \centering
    \subfloat[H-MG A]
    {
        \includegraphics[width = 0.95\columnwidth]{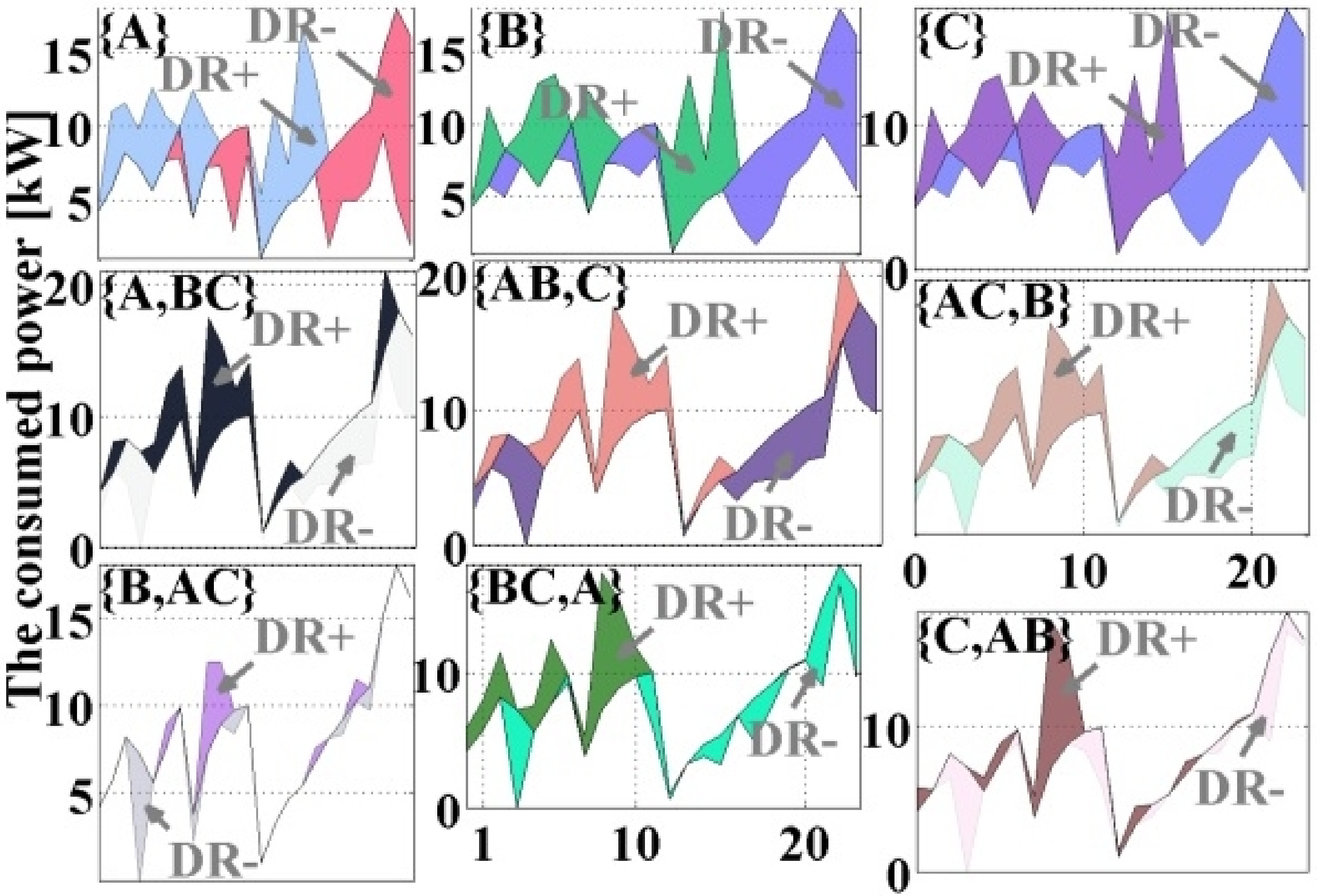}
        \label{fig:Fig_5_1}
    }\\
    \subfloat[H-MG B]
    {
        \includegraphics[width = 0.95\columnwidth]{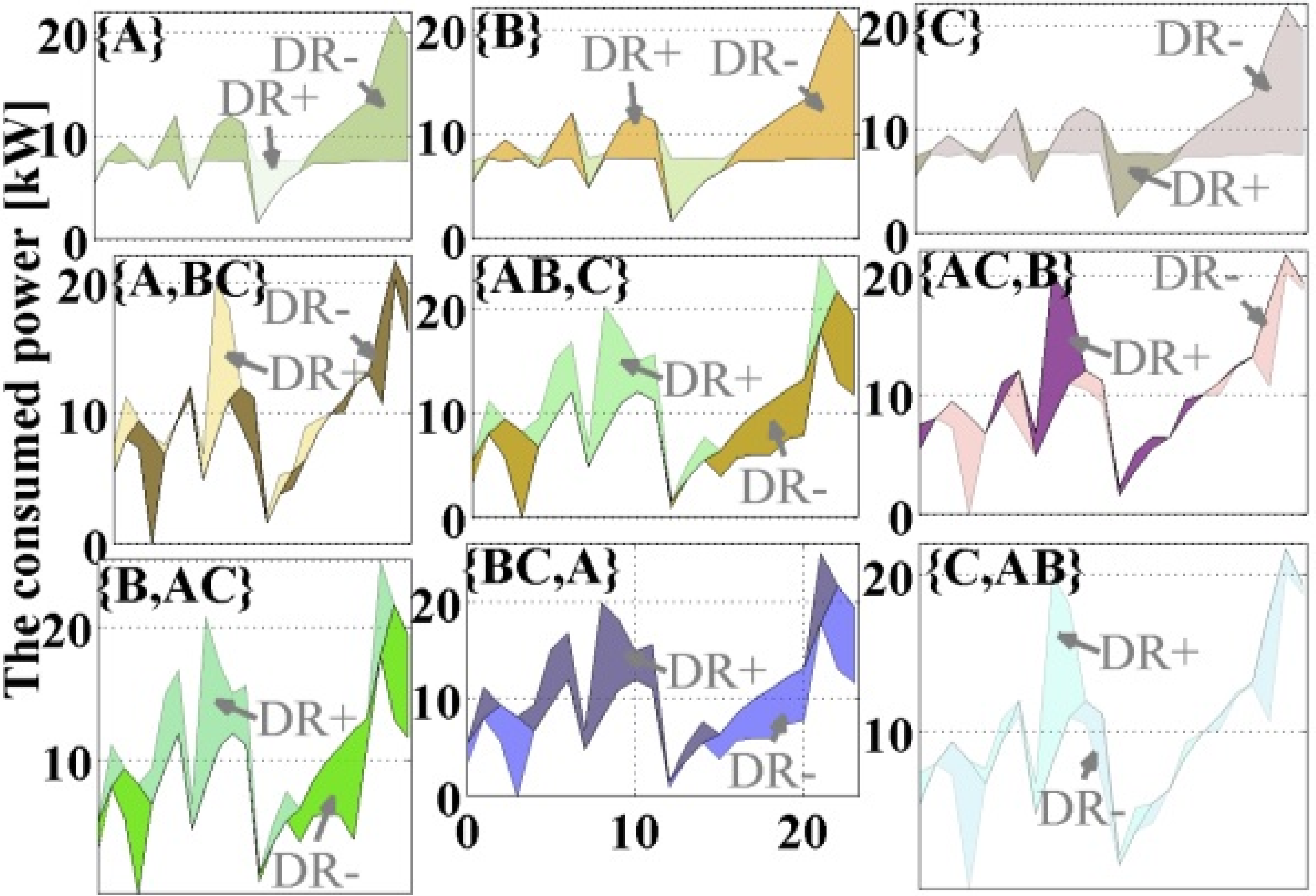}
        \label{fig:Fig_5_2}
    }\\
    \subfloat[H-MG C]
        {
            \includegraphics[width = 0.95\columnwidth]{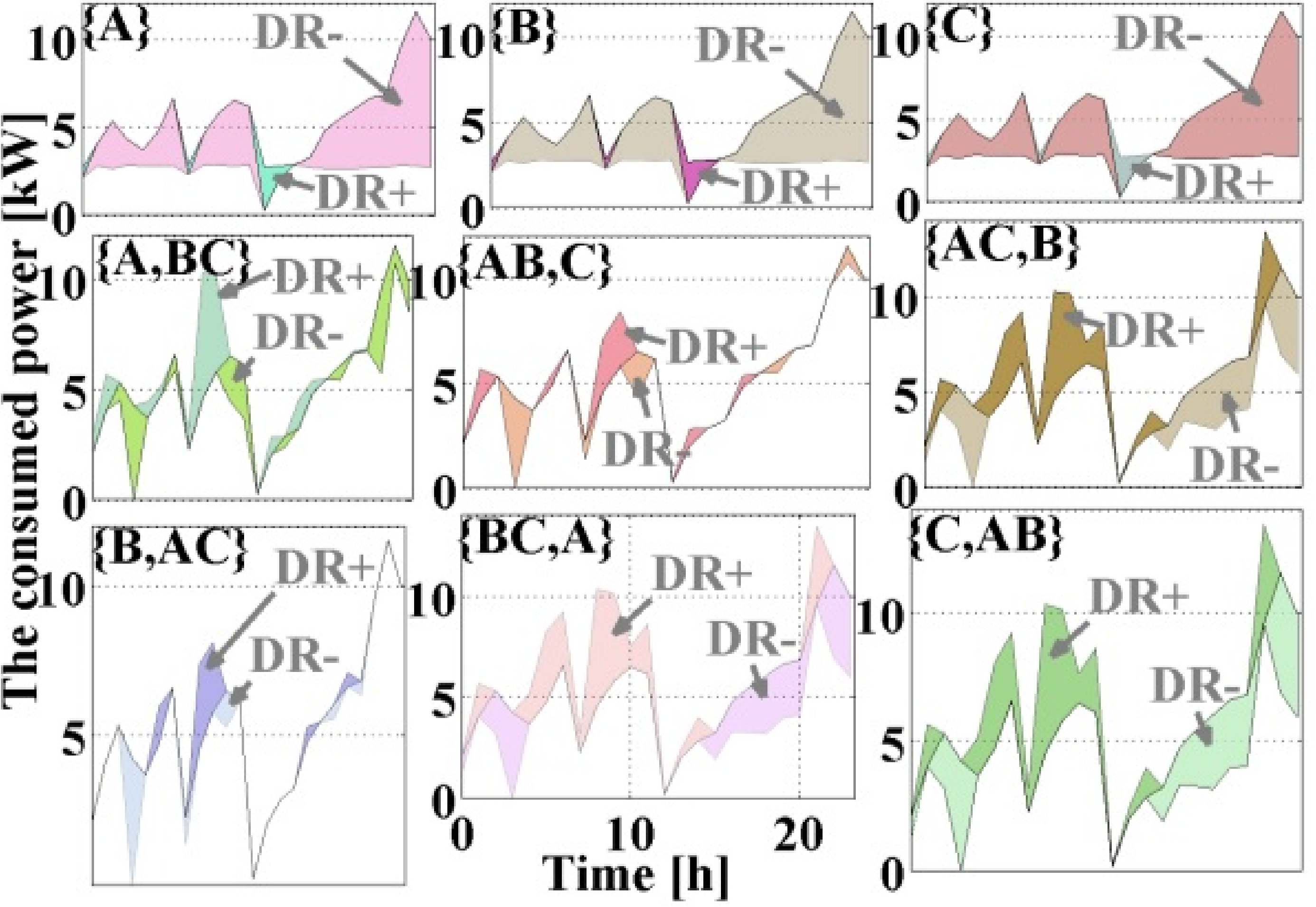}
            \label{fig:Fig_5_3}
        }
    \caption{DR+ and DR- in the H-MGs under different scenarios.}
    \label{fig:Load_Demand}
\end{figure}

In Figure.~\ref{fig:Load_Demand}, load profiles of H-MGs ({A}, {B} and {C}) are shown. For an independent operation of H-MGs, most of RLDs of H-MG A are shifted from the time intervals with higher MCP to the time intervals with lower MCP. The amount of load shifted forms a high share of the total load of H-MGs. More specifically, 55\% of the load is shifted from time intervals with higher MCP to other time intervals with lower MCP taking profit maximization for H-MGs owners. On the other hand, the energy consumption in such a figure has reduced significantly when H-MGs operate in coalitional structures. Such an energy consumption is the lowest (21\%) for a coalitional scenario of (\{B,AC\}). The energy consumption is at the lowest level (21\%) when the coalition scenario corresponds to (\{B, AC\}). In addition, the reduction in the amount of load shift is a result of a DSM program aimed at achieving a higher pay-off for consumers by considering (a) employing load shifting when the value of MCP is high, along with the maximum use of H-MG A interval resources, and (b) reducing the generation cost in an effective manner when load shifting is at a minimum. Moreover, the load profile of H-MGs in a coalition structure (\{A,BC\}) is the same as that of the alternative coalition structure (\{AB,C\}) and does not have a significant effect on the consumption level in H-MG A. This trend is completely different from the case in H-MG B. More specifically, during independent operation of H-MGs, the amount of load shift in H-MG B is at its least amount (almost 30\% of the total load during 24 hours). Therefore, forming a coalition among H-MGs would increase consumers's participation in the DR program that can reach almost 42\% to 50\%.\par

Such a reduction in the amount of load shift is a result of a DSM program  for reaching more pay-off for consumers by considering criteria such as load shifting when the value of MCP is high, the maximum use of H-MG A interval resources, and also the reduction in the value of generation cost in the best way and with the least amount of load shifting has taken place. Alternatively, the load profile of H-MGs in a coalitional structure (\{A,BC\}) is the same as that of such H-MGs in an alternative coalitional structure (\{AC,B\}), and does not have a significant effect on the consumption nature in H-MG A.

The least amount of load shifting is achieved when H-MGs B and A from a coalition in the lower-level of the BL-HIA structure, while having the objective function at an upper-level of the structure targeting maximizing the profit of H-MG C. Furthermore, these conditions are comparable for the \{AC, B\} coalition structure, having a similar nature. Under the previous conditions, a substantial share of the excess generation capacity is devoted to meeting H-MG C demand. As a result, a negligible part of such energy has been allocated for supplying responsive loads in H-MG B. It is important to clarify that in the case of H-MG C, the value of the total DR-, in H-MGs independent operating conditions, is significantly greater than the value corresponding to positive demand response (DR+) conditions. While only 17\% of time intervals, H-MG C had experienced a DR+ algorithm, such a figure would reach 83\% when a negative demand response (DR-) would be experienced. Such a trend in demand response (DR) is comparable to the scenario of coalition structures, where the amount of the total load shift, with the value of DR+ total load during daily performance, are close to each other in terms of the value. The consumer's participation percentage in the H-MG C has improved significantly by forming a coalition between H-MG B and A reaching more than 40\% of the times. Only in the coalitional structure \{B,AC\}, such a value can be minimum (21\%).

Furthermore, these conditions are also exactly similar for the \{AC,B\} coalitional structure and has a similar nature. Under the previous conditions, a big share of the amount of excess generation is spent supplying H-MG C demand.

It is important to clarify that for H-MG C, the value
of the total DR- in H-MGs independent operation conditions is much more than such a value when positive demand response (DR+) has taken place. 

Such a trend in demand response (DR) is quite similar to the scenario of a coalitional structures, where the amount of the total load shift and the value of DR+ total load during daily performance are close to each other in terms of the value.



\begin{figure}
\centering
\includegraphics[width= \columnwidth]{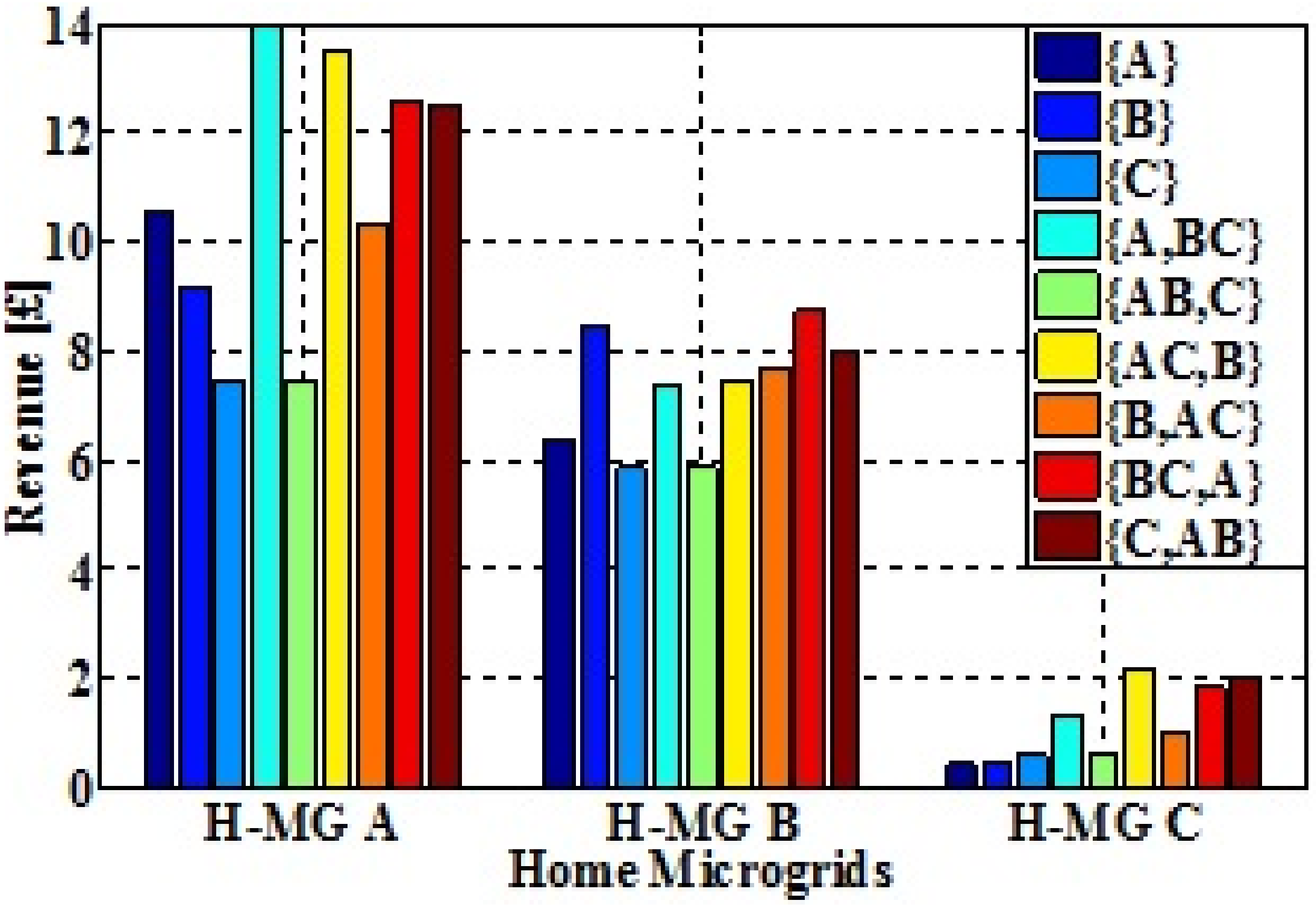}
\caption{H-MGs income under different scenarios.}
\label{fig:Fig_6}
\end{figure}

The increasing trend of each H-MG income during an independent and a coalitional performance with other H-MGs is shown in Figure.~\ref{fig:Fig_6}. As it is observed from this figure, each structure can be useful for one H-MG and meanwhile can have no benefit for other H-MGs. The best structure, which can be useful for H-MG A, results from forming a coalition among H-MG B and H-MG C excluding the participation of H-MG A in this coalition. These conditions can also  be useful for H-MG B on the condition of forming a coalition with H-MG C is in a higher priority of operation. For H-MG C, the highest income is experienced when this H-MG forms a coalition with H-MG A at an initial stage given that H-MG B works independently. Under these conditions, the income of H-MG A is close to the maximum value. For H-MG C, because of the lower the generated power, it is appropriate to form a coalition in all cases with other H-MGs. In all cases in which H-MG C has formed a coalition with other H-MGs, an increasing trend in the income is observed. In comparison, the income resulting from H-MG B when used independently, is significantly improved when compared to other configurations, such as coalition formation with other H-MGs. Furthermore, it is possible, in some cases, for coalition forming to have a detrimental effect on the H-MGs that form part of the coalition.  \par  
It is also observed that the coalition among H-MG A and H-MG B at the initial level leads to a significant reduction in the income independently obtained by this H-MG. Moreover, it is  desirable to prevent H-MG A from forming a coalition with H-MG B and negotiate with H-MG C to form the coalition.

In comparison, the income resulting from the independent performance of H-MG B is also significant compared to other cases (e.g., coalition formation with other H-MGs) and in some cases forming coalition is harmful for these H-MGs.\par

 For H-MG C, because of the lowness of the generated power, it is appropriate to form a coalition in all cases with other H-MGs.

\begin{figure}[ht!]
    \centering
    \subfloat[Electrical MCP]
    {
       \includegraphics[width = \columnwidth]{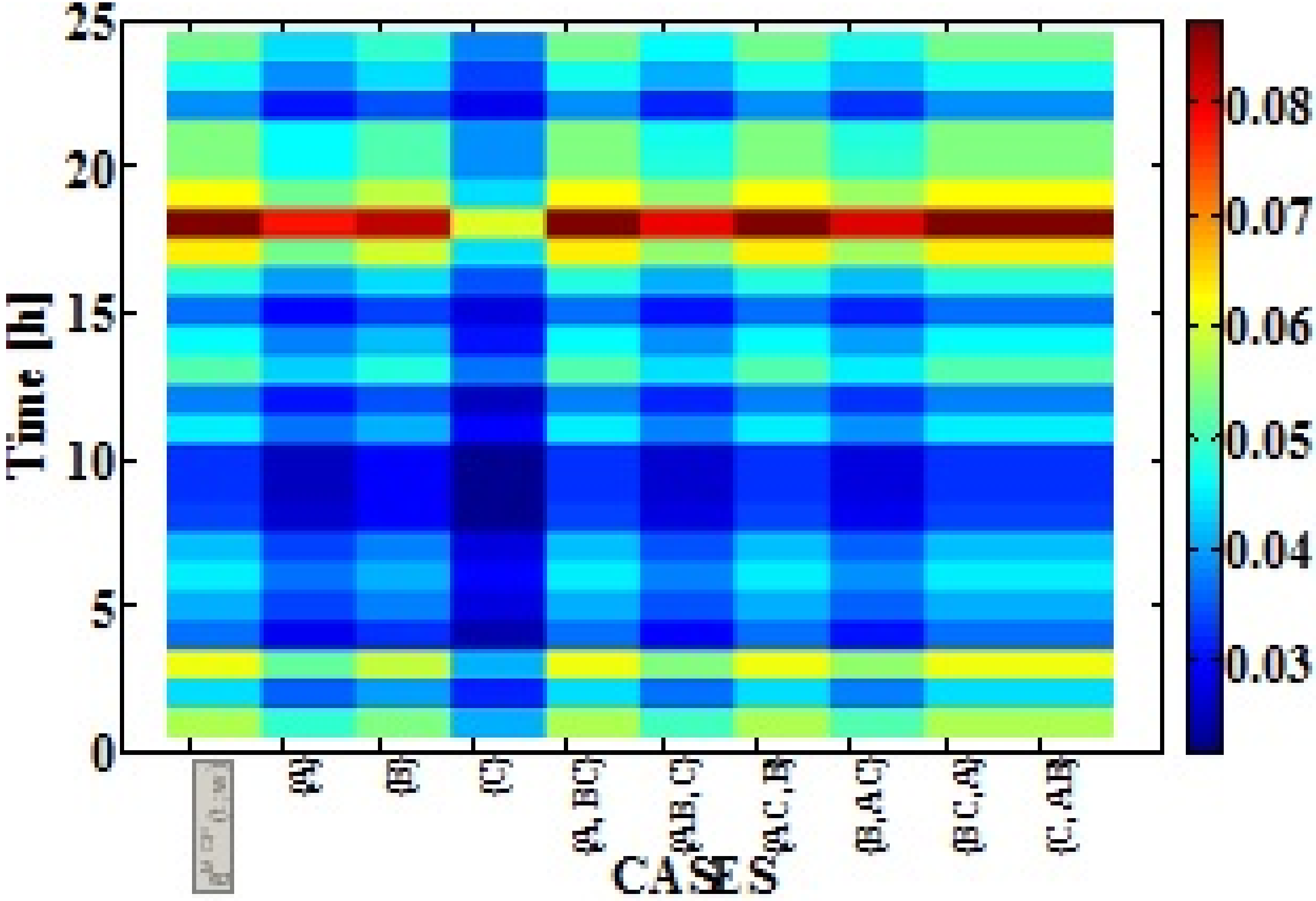}
        \label{fig:Fig_7_1}
    }\\
    \subfloat[Thermal MCP]
    {
        \includegraphics[width = \columnwidth]{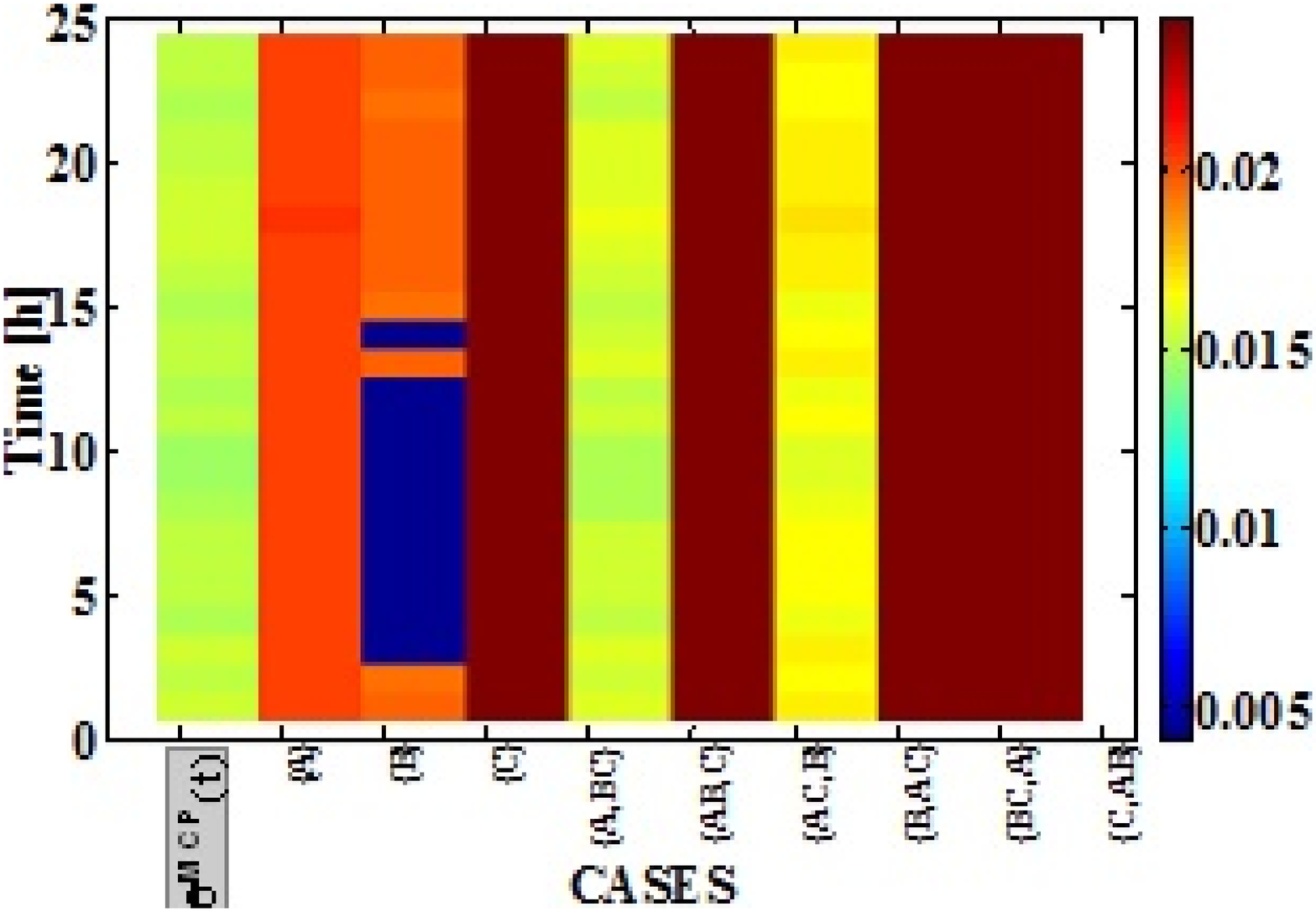}
        \label{fig:Fig_7_2}
    }\\
    \caption{The electrical and thermal MCP during the 24-hour system performance.}
    \label{fig:MCP}
\end{figure}

In Figure.~\ref{fig:Fig_7_1} and Figure.~\ref{fig:Fig_7_2} show the values of the electrical and thermal MCP, respectively. Although the average value of the electrical MCP in the case of an independent operation of H-MG C is at its minimum during the system's daily  performance, such values can be significantly improved when investigated at individual time intervals (i.e., one hour after forming a coalition among H-MGs). In some of time intervals, not only forming a coalition does not cause a degradation in the electrical MCP, but also a small increase in its value. Moreover, at certain intervals, its value may not change significantly when a coalition exists compared to the scenario where H-MGs work independently. The electrical MCP value in \{A,BC\} coalition is about 54\% of the times become more than its value in \{A,BC\} combination. That is why, no difference in values of the electrical MCP is observed for coalitions \{A,BC\}, \{AC,B\}, \{BC,A\} and \{C,AB\}. Furthermore, by changing the structure from \{A,BC\} to \{B,AC\}, about 33\% of the MCP value is reduced. \par
Such analysis also applies for the thermal MCP for the investigated structures. Finally, we can conclude from the simulation results that forming a coalition among H-MGs existing in one grid will not only have a significant effect on programming and regulating the value of the power generated by the generation resources but also can affect the change in the demand consumption and the behavior of consumers participating in the DR program with a cheaper MCP.


\section{Conclusion}

The paper presented an optimum development combined problem of the quantity in a deregulated electricity market environment. A methodology has been presented for investigating the possibility of increasing incomes of H-MGs, consumers and retailers existing in a multiple H-MGs. These Participants's performances had been properly modelled in the market environment. A H-MG programmer tries to increase its income as long as it is freely negotiating energy exchange with DER resources and its consumers. It can also put the possibility of forming a coalition with other H-MGs on its agenda. H-MGs seek to estimate the value of the power generated by DERs and also supply/demand bids to consumers. 
Meanwhile, the possibility of forming coalition among H-MGs with the aim of maximizing the income in an independent or a coalitional operation in a scheduling horizon is also investigated. In this way, the H-MGs encounter pool prices uncertainties, and the value of electrical and thermal loads. Furthermore, if the supply bid of one H-MG is not competitive enough, consumers may choose another H-MG for supplying their demand. For investigating how forming a coalition among H-MGs can affect the market behaviour and  H-MGs's gained income, different scenarios were presented. These scenarios were solved through a bilevel structure which can be transformed into a one NLP problem. The proposed model did not only present solutions of higher income achievements of each H-MG in an independent or a coalitional operation, but also provides the higher income/lower cost for each of the retailers/consumers relative to a single level model.\par 
The BL-HIA structure has presented an adequate framework for modelling both H-MGs reaction for a better participation in generation and effect on electricity price, and also competition increase among H-MGs and retailers. In the upper-level problem, H-MGs change their capacity with the aim of maximizing their income and by predicting the behaviour of other competitors (H-MGs) resulting  from the lower-level problem, and noting quantities and prices proposed by DERs and consumers. An Optimum pricing strategy was implemented to enable the market dynamic behavior on H-MGs decisions. Furthermore, a daily generation scheduling was presented. For a selected case study, an infinite number of Nash equilibrium was observed for the case where no players tend to unilaterally change their pricing strategies. In these obtained equilibrium points, the total expected profit of all players does not change. Although it is distributed among them.\par 
Simulation results showed that by forming a coalition among H-MGs, their profit, the demand value of supplied load and the DERs generated power existing in those H-MGs may change. Furthermore, computational simulations showed the convergence of the proposed model for solving real problems and simultaneously presenting solutions for raising H-MGs and retailers income and also reducing market clearing price. 
The following results can be extracted from the structure of the developed model:\par

\begin{itemize}
\item The bilevel model hierarchical structure for modelling the strategic behavior of each H-MG in reaction to the behavioral change and decision making of other H-MGs and their supply bid. Furthermore, the proposed structure can effectively encourage consumers to participate in the electricity market and affecting them using the DSM program.
\item It had been shown that energy exchange among H-MGs and retailers, in addition to increasing the profit of each player, would have a significant impact on levelling the load and reducing consumers's power consumption lack of support during the consumption peak period.
\end{itemize}

\textbf{Appendices}

\begin{figure}[ht!]
\centering
\includegraphics[width= \columnwidth]{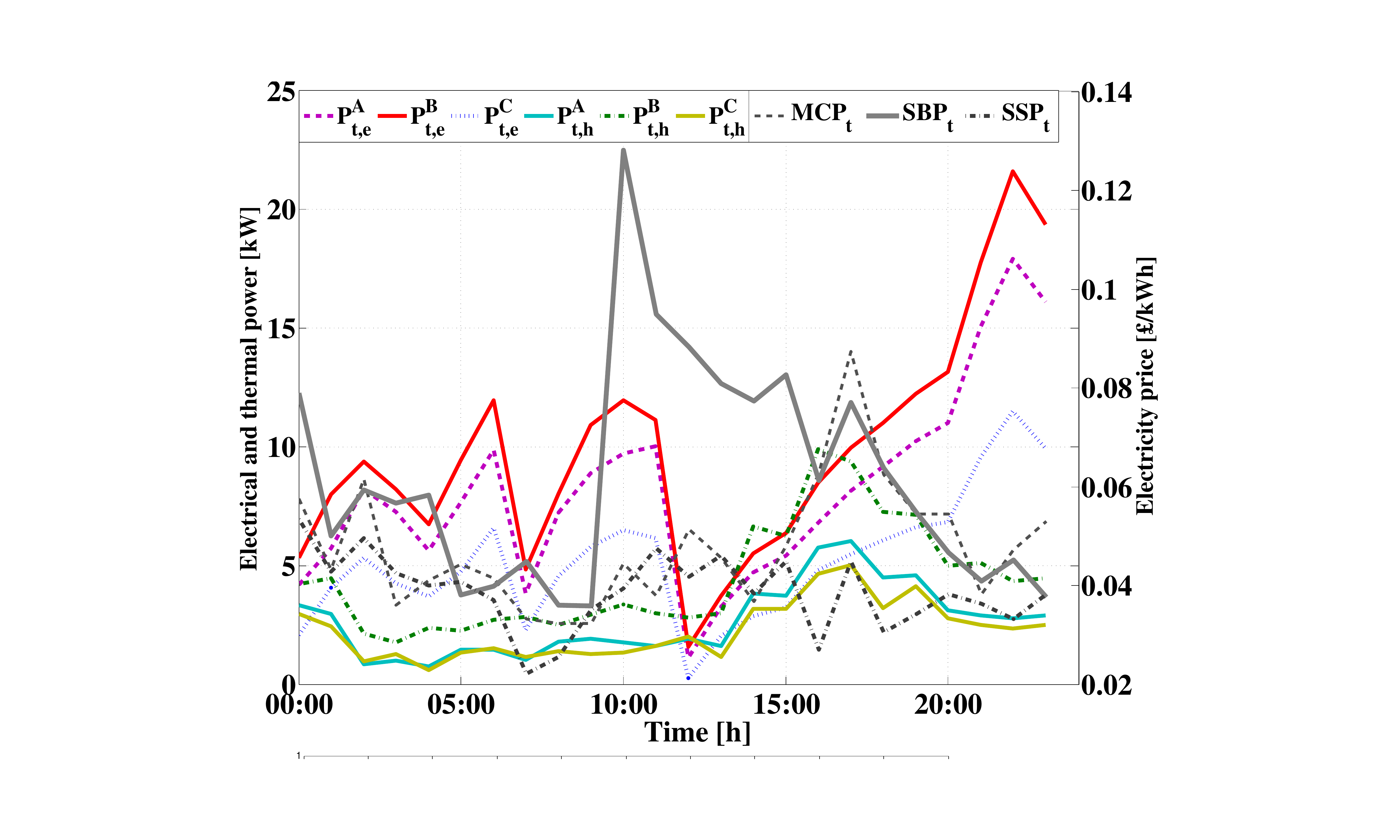}
\caption{The profiles of the electrical and thermal loads of the H-MGs and electricity prices for buying and selling.}
\label{fig:Load_pridect}
\end{figure}

\begin{figure}[ht!]
\centering
\includegraphics[width= \columnwidth]{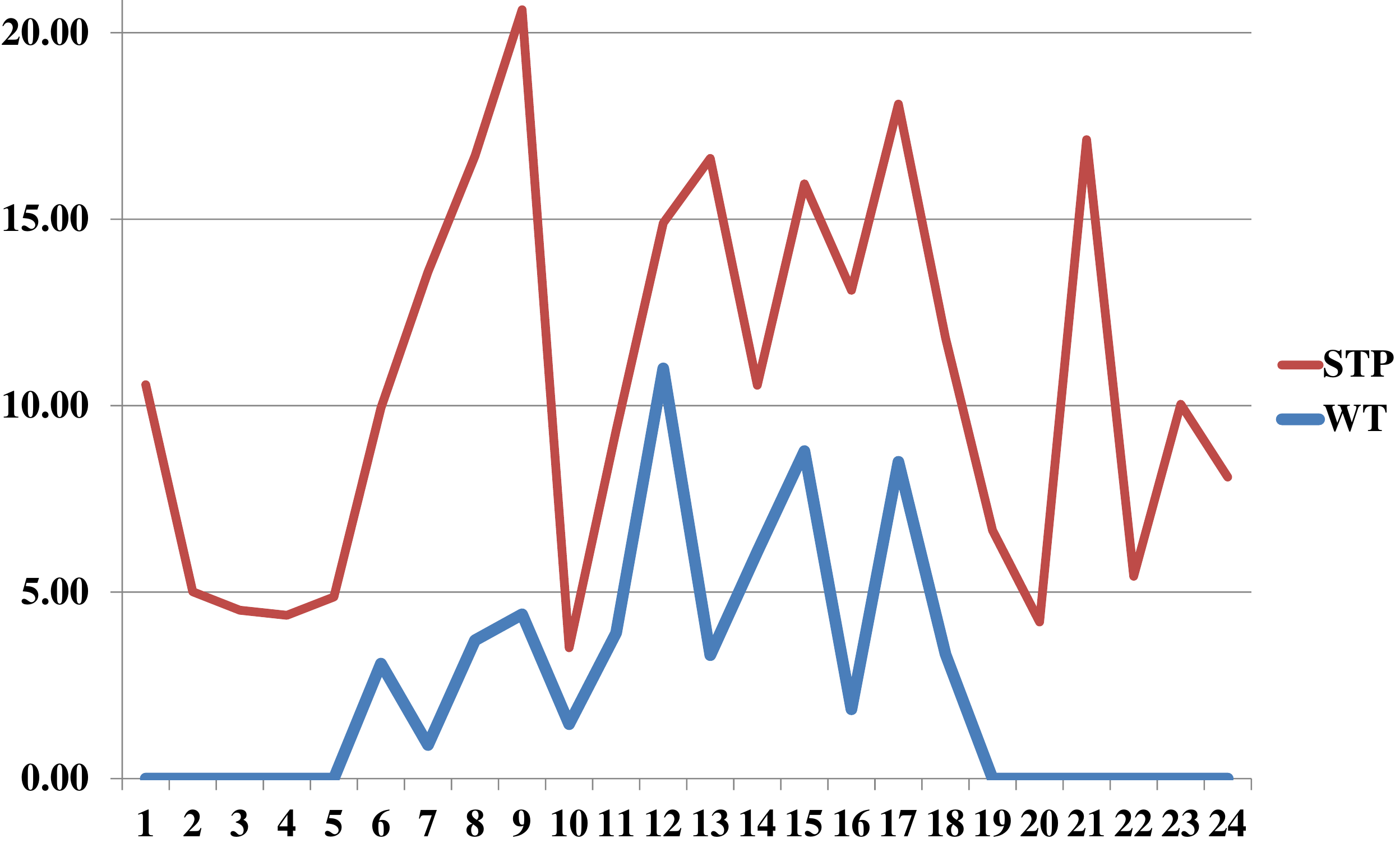}
\caption{Shape of output power wave generated by WT and STP.}
\label{fig:WTandSTP_Power}
\end{figure}

\textbf{Applying KKT conditions to the lower-level problem}

As stated before, since the lower-level problem is a linear problem, KKT conditions can be applied to it. The constraints obtained from the derivative of Lagrange expression relative to lower-level variables include the following relations:

\noindent
\begin{equation}
\label{eq:67}
\begin{array}{l}
\frac{\partial L}{\partial P_{\text{itt'w}}^{\prime D,e}}
= 0.5
\displaystyle{\sum_{t'=1}^{T}} 
\left(
\tilde{\lambda}_{\text{tw}}^{\text{MCP,e}}
- \tilde{\lambda}_{\text{t'w}}^{\text{MCP,e}}
\right)
- \gamma_{\text{itw}^{D,e}}|_{P_{itt'w}^{\prime D,e}>0} \\
\overline{\eta}_{\text{itw}}^{D,e}|_{P_{itt'w}^{\prime D,e}>0}
+ \eta_{\text{itt'w}}^{\prime D,e}
+ \eta_{\text{it'tw}}^{\prime D,e}= 0:
\;\;
\eta_{\text{itt'w}}^{D,e\min}
\end{array}
\end{equation}


\noindent
\begin{equation}
\label{eq:68}
\frac{\partial L}{\partial P_{\text{itw}}^{D,e}}
= \gamma_{\text{itw}}^{D,e}
+ \lambda_{\text{tw}}^{h}
= 0
\end{equation}

\noindent
\begin{equation}
\label{eq:69}
\frac{\partial L}{\partial P_{\text{ijtw}}^{\text{CHP,e}}}
= \pi_{\text{ijtw}}^{\text{CHP,e}}
+ \overline{\eta}_{\text{ijtw}}^{\text{CHP,e}}
- \underline{\eta}_{\text{ijtw}}^{\text{CHP,e}}
+ \gamma_{\text{ijtw}}^{\text{CHP,e}}
- \lambda_{\text{tw}}^{e}
= 0
\end{equation}

\noindent
\begin{equation}
\label{eq:70}
\frac{\partial L}{\partial P_{\text{ijtw}}^{\text{CHP,h}}}
= \pi_{\text{ijtw}}^{\text{CHP,h}}
- \frac{\pi^{\text{FU}}}{N_{i}^{\text{CHP,h}}}
+ \overline{\eta}_{\text{ijtw}}^{\text{CHP,h}}
- \underline{\eta}_{\text{ijtw}}^{\text{CHP,h}}
- \frac{\gamma_{\text{ijtw}}^{\text{CHP,e}}\zeta_{\text{ij}}^{\text{CHP,e}}}{\zeta_{\text{ij}}^{\text{CHP,h}}}
- \lambda{\text{tw}}^{h}= 0
\end{equation}

\noindent
\begin{equation}
\label{eq:71}
\frac{\partial L}{\partial P_{\text{ijtw}}^{\text{WT}}}
= \pi_{\text{ijtw}}^{\text{WT}}
+ \overline{\gamma}_{\text{tw}}^{\text{WT,e}}
- \underline{\gamma}_{\text{tw}}^{\text{WT,e}}
- \lambda_{\text{tw}}^{e}= 0
\end{equation}

\noindent
\begin{equation}
\label{eq:72}
\frac{\partial L}{\partial P_{\text{ijtw}}^{\text{ES}}}
= \pi_{\text{ijtw}}^{\text{ES}}
+ \overline{\eta}_{\text{tw}}^{\text{ES,e}}
- \underline{\eta}_{\text{tw}}^{\text{ES,e}}
- \frac{\gamma_{\text{ijtw}}^{\text{ES,e}}}{\overline{P}_{ij}^{\text{ES}}}
- \lambda_{\text{tw}}^{e}
= 0
\end{equation}

\noindent
\begin{equation}
\label{eq:73}
\frac{\partial L}{\partial \text{SOC}_{\text{ijtw}}^{\text{ES}}}
= \gamma_{\text{ijtw}}^{\text{ES,e}}
- \gamma_{ij(t-1)w}^{\text{ES,e}}
+ \overline{\eta}_{\text{ijtw}}^{\text{ES,SOC}}
- \underline{\eta}_{\text{ijtw}}^{\text{ES,SOC}}
+ \gamma_{\text{ijw}}^{\prime \text{ES}}|_{t=1}
+ \gamma_{\text{ijw}}^{\prime \prime \text{ES}}|_{t=24}
= 0
\end{equation}

\noindent
\begin{equation}
\label{eq:74}
\frac{\partial L}{\partial P_{\text{ijtw}}^{\text{EB,h}}}
+ \gamma_{\text{tw}}^{\text{EB,h}}
+ \overline{\eta}_{\text{ijtw}}^{\text{EB,h}}
- \underline{\eta}_{\text{ijtw}}^{\text{EB,h}}
- \lambda_{\text{tw}}^{h}
= 0
\end{equation}

\noindent
\begin{equation}
\label{eq:75}
\frac{\partial L}{\partial P_{\text{ijtw}}^{\text{EB,e}}}
= - \pi_{\text{ijtw}}^{\text{EB,e}}
- \zeta_{\text{ij}}^{\text{EB}}
\times
\gamma_{\text{tw}}^{\text{EB,h}}
+ \lambda_{\text{tw}}^{e}
= 0
\end{equation}

\noindent
\begin{equation}
\label{eq:76}
\frac{\partial L}{\partial P_{\text{ijtw}}^{\text{EHP,h}}}
= \pi_{\text{ijtw}}^{\text{EHP,h}}
+ \gamma_{\text{ijtw}}^{\text{EHP,h}}
+ \overline{\eta}_{\text{ijtw}}^{\text{EHP,h}}
- \underline{\eta}_{\text{ijtw}}^{\text{EHP,h}}
- \lambda_{\text{tw}}^{h}
= 0
\end{equation}

\noindent
\begin{equation}
\label{eq:77}
\frac{\partial L}{\partial P_{\text{ijtw}}^{\text{EHP,e}}}
= -\pi_{\text{ijtw}}^{\text{EHP,e}}
- \text{COP}^{\text{EHP}}
\times
\gamma_{\text{ijtw}}^{\text{EHP,h}}
+ \lambda_{\text{tw}}^{e}
= 0
\end{equation}

\noindent
\begin{equation}
\label{eq:78}
\frac{\partial L}{\partial P_{\text{ijtw}}^{\text{STP}}}
= \pi_{\text{ijtw}}^{\text{STP}}
+ \overline{\eta}_{\text{ijtw}}^{\text{STP,h}}
- \underline{\eta}_{\text{ijtw}}^{\text{STP,h}}
- \lambda_{\text{tw}}^{h}
= 0
\end{equation}

\noindent
\begin{equation}
\label{eq:79}
\frac{\partial L}{\partial P_{\text{ijtw}}^{\text{GB}}}
= \pi{\text{ijtw}}^{\text{GB}}
- \frac{\pi^{\text{FU}}}{N_{\text{ij}}^{\text{GB}}}
+ \overline{\eta}_{\text{ijtw}}^{\text{GB,h}}
- \underline{\eta}_{\text{ijtw}}^{\text{GB,h}}
- \lambda_{\text{tw}}^{h}
= 0
\end{equation}

\noindent
\begin{equation}
\label{eq:80}
\frac{\partial L}{\partial P_{\text{ijtw}}^{\text{TES}}}
= \pi_{\text{ijtw}}^{\text{TES}}
+ \overline{\eta}_{\text{ijtw}}^{\text{TES,h}}
- \underline{\eta}_{\text{ijtw}}^{\text{TES,h}}
- \frac{\gamma_{\text{ijtw}}^{\text{TES,h}}}{\overline{P}_{\text{ij}}^{\text{TES}}}
- \lambda_{\text{tw}}^{h}
= 0
\end{equation}

\noindent
\begin{equation}
\label{eq:81}
\begin{array}{l}
\frac{\partial L}{\partial \text{SOC}_{\text{ijtw}}^{\text{TES}}}
= \gamma_{\text{ijtw}}^{\text{TES}}
- \gamma_{\text{ij(t-1)w}}^{\text{TES,h}}
\times
\zeta_{\text{ij}}^{\text{TES}}
+ \overline{\eta}_{\text{ijtw}}^{\text{TES,SOC}},
\\
\underline{\eta}_{\text{ijtw}}^{\text{TES,SOC}}
+ \gamma_{\text{ijw}}^{\prime \text{TES}}|_{t=1}
+ \gamma_{\text{ijw}}^{\prime \prime \text{TES}}|_{t=24}
= 0
\end{array}
\end{equation}

The complementary constraints obtained from the inequality constraints of lower-level problem are in the form of following relations:

\noindent
\begin{equation}
\label{eq:82}
0
\leq
\left(
\overline{P}_{\text{ij}}^{\text{CHP,e}}
- P_{\text{ijtw}}^{\text{CHP,e}}
\right)
\perp
\overline{\eta}_{\text{ijtw}}^{\text{CHP,e}}
\geq
0
\end{equation}

\noindent
\begin{equation}
\label{eq:83}
0
\leq
\left(
P_{\text{ijtw}}^{\text{CHP,e}}
- 
\underline{P}_{\text{ij}}^{\text{CHP,e}}
\right)
\perp
\underline{\eta}_{\text{ijtw}}^{\text{CHP,e}}
\geq
0
\end{equation}

\noindent
\begin{equation}
\label{eq:84}
0
\leq
\left(
\overline{P}_{\text{ij}}^{\text{CHP,h}}
- P_{\text{ijtw}}^{\text{CHP,h}}
\right)
\perp
\overline{\eta}_{\text{ijtw}}^{\text{CHP,h}}
\geq
0
\end{equation}

\noindent
\begin{equation}
\label{eq:85}
0
\leq
\left(
P_{\text{ijtw}}^{\text{CHP,h}}
- 
\underline{P}_{\text{ij}}^{\text{CHP,h}}
\right)
\perp
\underline{\eta}_{\text{ijtw}}^{\text{CHP,h}}
\geq
0
\end{equation}

\noindent
\begin{equation}
\label{eq:86}
0
\leq
P_{\text{ijtw}}^{\text{WT}}
\perp
\underline{\gamma}_{\text{tw}}^{\text{WT,e}}
\geq
0
\end{equation}

\noindent
\begin{equation}
\label{eq:87}
0
\leq
\left(
\overline{P}_{\text{ij}}^{\text{WT}}
- P_{\text{ijtw}}^{\text{WT}}
\right)
\perp
\overline{\gamma}_{\text{ijtw}}^{\text{WT,e}}
\geq
0
\end{equation}

\noindent
\begin{equation}
\label{eq:88}
0
\leq
\left(
P_{\text{ijtw}}^{\text{ES}}
- \underline{P}_{\text{ij}}^{\text{ES}}
\right)
\perp
\underline{\eta}_{\text{ijtw}}^{\text{ES,e}}
\geq
0
\end{equation}

\noindent
\begin{equation}
\label{eq:89}
0
\leq
\left(
\overline{P}_{\text{ij}}^{\text{ES}}
- P_{\text{ijtw}}^{\text{ES}}
\right)
\perp
\overline{\eta}_{\text{ijtw}}^{\text{ES,e}}
\geq
0
\end{equation}

\noindent
\begin{equation}
\label{eq:90}
0
\leq
\left(
\text{SOC}_{\text{ijtw}}^{\text{ES}}
- \underline{P}_{\text{ij}}^{\text{ES}}
\right)
\perp
\underline{\eta}_{\text{ijtw}}^{\text{ES,SOC}}
\geq
0
\end{equation}

\noindent
\begin{equation}
\label{eq:91}
0
\leq
\left(
\overline{\text{SOC}}_{\text{ij}}^{\text{ES}}
- \text{SOC}_{\text{ijtw}}^{\text{ES}}
\right)
\perp
\underline{\eta}_{\text{ijtw}}^{\text{ES,SOC}}
\geq
0
\end{equation}

\noindent
\begin{equation}
\label{eq:92}
0
\leq
P_{\text{ijtw}}^{\text{EB,h}}
\perp
\underline{\eta}_{\text{ijtw}}^{\text{EB,h}}
\geq
0
\end{equation}

\noindent
\begin{equation}
\label{eq:93}
0
\leq
\left(
\overline{P}_{\text{ij}}^{\text{EB,h}}
- P_{\text{ijtw}}^{\text{EB,h}}
\right)
\perp
\overline{\eta}_{\text{ijtw}}^{\text{EB,h}}
\geq
0
\end{equation}

\noindent
\begin{equation}
\label{eq:94}
0
\leq
P_{\text{ijtw}}^{\text{EHP,h}}
\perp
\underline{\eta}_{\text{ijtw}}^{\text{EHP,h}}
\geq
0
\end{equation}

\noindent
\begin{equation}
\label{eq:95}
0
\leq
\left(
\overline{P}_{\text{ij}}^{\text{EHP,h}}
- P_{\text{ijtw}}^{\text{EHP,h}}
\right)
\perp
\overline{\eta}_{\text{ijtw}}^{\text{EHP,h}}
\geq
0
\end{equation}

\noindent
\begin{equation}
\label{eq:96}
0
\leq
\left(
P_{\text{ijtw}}^{\text{TES}}
- \underline{P}_{\text{ij}}^{\text{TES}}
\right)
\perp
\underline{\eta}_{\text{ijtw}}^{\text{TES,h}}
\geq
0
\end{equation}

\noindent
\begin{equation}
\label{eq:97}
0
\leq
\left(
\overline{P}_{\text{ij}}^{\text{TES}}
- P_{\text{ijtw}}^{\text{TES}}
\right)
\perp
\overline{\eta}_{\text{ijtw}}^{\text{TES,h}}
\geq
0
\end{equation}

\noindent
\begin{equation}
\label{eq:98}
0
\leq
\left(
\text{SOC}_{\text{ijtw}}^{\text{TES}}
- \underline{P}_{\text{ij}}^{\text{TES}}
\right)
\perp
\underline{\eta}_{\text{ijtw}}^{\text{TES,SOC}}
\geq
0
\end{equation}

\noindent
\begin{equation}
\label{eq:99}
0
\leq
\left(
\overline{\text{SOC}}_{\text{ij}}^{\text{TES}}
- \text{SOC}_{\text{ijtw}}^{\text{TES}}
\right)
\perp
\underline{\eta}_{\text{ijtw}}^{\text{TES,SOC}}
\geq
0
\end{equation}

\noindent
\begin{equation}
\label{eq:100}
0
\leq
P_{\text{ijtw}}^{\text{GB}}
\perp
\underline{\eta}_{\text{ijtw}}^{\text{GB,h}}
\geq
0
\end{equation}

\noindent
\begin{equation}
\label{eq:101}
0
\leq
\left(
\overline{P}_{\text{ij}}^{\text{GB}}
- P_{\text{ijtw}}^{\text{GB}}
\right)
\perp
\overline{\eta}_{\text{ijtw}}^{\text{GB,h}}
\geq
0
\end{equation}

\noindent
\begin{equation}
\label{eq:102}
0
\leq
P_{\text{ijtw}}^{\text{STP}}
\perp
\underline{\eta}_{\text{ijtw}}^{\text{STP,h}}
\geq
0
\end{equation}

\noindent
\begin{equation}
\label{eq:103}
0
\leq
\left(
\overline{P}_{\text{ij}}^{\text{STP}}
- P_{\text{ijtw}}^{\text{STP}}
\right)
\perp
\overline{\eta}_{\text{ijtw}}^{\text{STP,h}}
\geq
0
\end{equation}

\noindent
\begin{equation}
\label{eq:104}
0
\leq
P_{\text{iktw}}^{\text{e+}}
\perp
\underline{\eta}_{\text{iktw}}^{\text{Grid,e}}
\geq
0
\end{equation}

\noindent
\begin{equation}
\label{eq:105}
0
\leq
\left(
\xi
\times
\left(
P_{\text{ijtw}}^{\text{CHP,e}}
+ P_{\text{ijtw}}^{\text{WT}}
+ P_{\text{ijtw}}^{\text{ES}}
- P_{\text{iktw}}^{\text{e+}}
\right)
\right)
\perp
\overline{\eta}_{\text{iktw}}^{Grid,e}
\geq
0
\end{equation}




\section*{References}
\bibliographystyle{model3-num-names}
\bibliography{neighbourhood} 

\cleardoublepage
\end{document}